\documentclass[11pt]{article}

\usepackage[usenames,dvipsnames]{pstricks}
\usepackage{epsfig}
\usepackage{pst-grad} 
\usepackage[space]{grffile} 
\usepackage{etoolbox} 
\makeatletter 
\patchcmd\Gread@eps{\@inputcheck#1 }{\@inputcheck"#1"\relax}{}{}
\makeatother

\usepackage{graphicx}
\usepackage{amsmath}
\usepackage{psfrag}
\usepackage{color,soul}
\usepackage{natbib}
\usepackage{siunitx}
\usepackage{float}
\usepackage{overpic}
\usepackage{tikz}
\usepackage{multirow}
\usepackage{bm}
\usepackage{textgreek}
\usepackage{url}
\usepackage{amssymb}
\usepackage{wasysym}
\usepackage{xfrac}
\usepackage{setspace}
\usepackage{gensymb}
\usepackage{mathtools}
\usepackage{arydshln}
\usepackage{algorithm}
\usepackage{algpseudocode}

\usepackage[noblocks]{authblk}
\usepackage[margin=1in]{geometry}
\usepackage[title]{appendix}

\newcommand\affiliation[1]{\gdef\@affiliation{\let\aff\aff@inst#1}}
\gdef\@affiliation{}
\def\email#1{Email address for correspondence: #1}
\def\aff#1{\ignorespaces\textsuperscript{#1}}
\def\corresp#1{\unskip\thanks{#1}}
\numberwithin{equation}{section}

\renewenvironment{abstract}
{\begin{quote}
\noindent \rule{\linewidth}{.5pt}\par{\bfseries \abstractname.}}
{\medskip\noindent \rule{\linewidth}{.5pt}
\end{quote}
}

  \DeclareTextFontCommand\textsfi{\usefont{OT1}{cmss}{m}{sl}}
  \DeclareMathAlphabet\mathsfi            {OT1}{cmss}{m}{sl}
  \DeclareTextFontCommand\textsfb{\usefont{OT1}{cmss}{bx}{n}}
  \DeclareMathAlphabet\mathsfb            {OT1}{cmss}{bx}{n}
  \DeclareTextFontCommand\textsfbi{\usefont{OT1}{cmss}{m}{sl}}
  \DeclareMathAlphabet\mathsfbi            {OT1}{cmss}{m}{sl}
  \IfFileExists{t1phv.fd}
  {\DeclareTextFontCommand\textsfbi{\usefont{T1}{phv}{b}{it}}
  \DeclareMathAlphabet\mathsfbi            {T1}{phv}{b}{it}
  }{}
  \IfFileExists{ot1phv.fd}
  {\DeclareTextFontCommand\textsfbi{\usefont{OT1}{phv}{b}{it}}
  \DeclareMathAlphabet\mathsfbi            {OT1}{phv}{b}{it}
  }{}

\newcommand{\ket}[1]{\boldsymbol{#1}}  
\newcommand{\tc}[1]{\mathsfbi{#1}}	    
\newcommand{\td}[1]{ {\bf #1} }			
\newcommand{\vc}[1]{\boldsymbol{#1}}    
\newcommand{\vd}[1]{ {\bf #1} }			
\newcommand{\vdh}[1]{ \hat{\ket{#1}}}
\newcommand{\citer}[1]{Ref. \citep{#1}}    
\newcommand{\citers}[1]{Refs. \citep{#1}}

\makeatletter
\newcommand{\lambdabar}{{\mathchoice
  {\smash@bar\textfont\displaystyle{0.25}{1.2}\lambda}
  {\smash@bar\textfont\textstyle{0.25}{1.2}\lambda}
  {\smash@bar\scriptfont\scriptstyle{0.25}{1.2}\lambda}
  {\smash@bar\scriptscriptfont\scriptscriptstyle{0.25}{1.2}\lambda}
}}
\newcommand{\smash@bar}[4]{%
  \smash{\rlap{\raisebox{-#3\fontdimen5#10}{$\m@th#2\mkern#4mu\mathchar'26$}}}%
}
\makeatother

\DeclareSymbolFont{matha}{OML}{txmi}{m}{it}
\DeclareMathOperator*{\argmax}{arg\,max}
\DeclareMathOperator*{\fint}{\frac{1}{T}\int_{0}^{T}}

\usepackage{hyperref}

\definecolor{darkblue}{rgb}{0,0,0.80}

\hypersetup{
    colorlinks=true,
    linkcolor={darkblue},
    citecolor={darkblue},
    filecolor=black,
    urlcolor={darkblue},
    plainpages=false,
}

\providecommand{\keywords}[1]
{

  \textbf{Key words: } #1
}

\title{\bf Linear model reduction using spectral proper orthogonal decomposition}

\author[1]{\bf Peter Frame\corresp{\email{pframe@umich.edu}}}
\author[2]{\bf Cong Lin}
\author[2]{\bf Oliver Schmidt}
\author[1]{\bf Aaron Towne}

\affil[1]{\normalsize Department of Mechanical Engineering, University of Michigan, Ann Arbor, MI, USA  }
\affil[2]{\normalsize Department of Mechanical and Aerospace Engineering, University of California, San Diego, CA, USA}

\date{}
\begin{document}
\maketitle

\begin{abstract}
Most model reduction methods reduce the state dimension and then temporally evolve a set of coefficients that encode the state in the reduced representation. In this paper, we instead employ an efficient representation of the entire trajectory of the state over some time interval of interest and then solve for the static coefficients that encode the trajectory on the interval. We use spectral proper orthogonal decomposition (SPOD) modes, which are provably optimal for representing long trajectories and substantially outperform any representation of the trajectory in a purely spatial basis (e.g., POD). We develop a method to solve for the SPOD coefficients that encode the trajectories for forced linear dynamical systems given the forcing and initial condition, thereby obtaining the accurate prediction of the dynamics afforded by the SPOD representation of the trajectory. The method, which we refer to as spectral solution operator projection (SSOP), is derived by projecting the general time-domain solution for a linear time-invariant system onto the SPOD modes. We demonstrate the new method using two examples: a linearized Ginzburg-Landau equation and an advection-diffusion problem. In both cases, the error of the proposed method is orders of magnitude lower than that of POD-Galerkin projection and balanced truncation. The method is also fast, with CPU time comparable to or lower than both benchmarks in our examples. Finally, we describe a data-free space-time method that is a derivative of the proposed method and show that it is also more accurate than balanced truncation in most cases.  \\  
\end{abstract}


\section{Introduction}
\label{Sec:Intro}

The expense of many modern computational models can prohibit their use in applications where speed is required. In a design optimization problem, for example, many simulations with different boundary conditions or parameters must be performed. In control applications, simulations may need to be conducted in real time to inform actuation. Model reduction techniques strive to deliver the orders-of-magnitude speedup necessary to enable adequately fast simulation for these and other problems with only a mild sacrifice in accuracy. 

The great majority of model reduction methods employ the following two-step strategy: (i) find an accurate compression of the state of the system at a particular time and (ii) find equations that evolve the coefficients that represent the state in this representation. The POD-Galerkin method \citep{Aubry88,Noack03,Rowley04}, perhaps the most widely used starting point for model reduction, is representative of this approach. The proper orthogonal decomposition (POD) modes are an efficient means of representing the state in that with relatively few POD coefficients, the state can often be represented to high accuracy. In a POD-Galerkin reduced-order model (ROM), these coefficients are then evolved by projecting the governing equations into the space of POD modes, yielding a much smaller dynamical system to evolve. Many alternative choices have been explored for both steps. Examples of the compression step include using balanced truncation modes \citep{Moore81} and autoencoders \citep{Lee20,gonzalez18}, and examples of deriving the equations in the reduced space include using Petrov-Galerkin projections \citep{Carlberg10,Carlberg13,Otto22} and learning the equations from data \cite{Peherstorfer16,Padovan24}. All of these approaches to model reduction, however, fall within the two-step strategy outlined above.
\\

We have investigated a different approach in this work: instead of representing the state (at a particular time) in a reduced manner, we instead employ a reduced representation for the entire trajectory, i.e., the state's evolution for some time interval. Whereas POD modes are the most efficient (linear) representation of the state, they are far from the most efficient representation of trajectories. This is true intuitively -- to represent a trajectory with POD modes, one has to specify the POD coefficients for each time step along the trajectory, but from one time step to the next, the POD coefficients are highly correlated. The analog of POD for entire trajectories is space-time POD \citep{Lumley67,Schmidt_Schmid19,Frame23}. Space-time POD modes are themselves time- and space-dependent, so to represent a trajectory, they are weighted by static coefficients. These modes are the most efficient linear representation of trajectories in the sense that to represent a trajectory to some desired accuracy, fewer degrees of freedom are needed if the trajectory is represented with space-time POD modes than any other linear encoding scheme. Unfortunately, space-time POD modes have a number of characteristics that make them undesirable for model reduction; computing them requires much training data, storing them is memory intensive, and computing space-time inner products, which would be necessary in a space-time ROM method, is expensive.
\\

Fortunately, an efficient space-time basis that does not share the undesirable properties of space-time POD modes exists. Spectral POD (SPOD) modes are most naturally formulated as a POD in the frequency domain. More precisely, at every temporal frequency, there exists a set of spatial modes that optimally represent the spatial structure at that frequency. These modes are the SPOD modes, and they may be thought of as space-time modes where each spatial mode $\boldsymbol{\psi}_{k,j}$ at frequency $\omega_k$ has the time dependence $e^{i\omega_k t}$. Each mode is associated with an energy, and these energies may be compared across frequencies; for example, the second mode at one frequency may have more energy than the first mode at another frequency. The fact that motivates this work is that the most energetic SPOD modes are also an excellent basis for representing trajectories. In fact, SPOD modes converge to space-time POD modes as the time interval becomes long, so for long intervals, the representation of a trajectory with SPOD modes is nearly as accurate (on average) as the space-time POD representation, which is optimal among all linear representations \citep{Frame23}.
\\

With this motivation, the goal of this work is to develop an algorithm to solve quickly for the SPOD coefficients that represent a trajectory for forced linear dynamical systems given the initial condition and forcing. If these coefficients can be obtained accurately, then the resulting error will be substantially lower than that of a POD-Galerkin model with the same number of modes. The method works as follows. The SPOD coefficients at a given frequency are related to the (temporal) Fourier transform of the state at the same frequency, which in turn is related to the forcing and initial condition. We derive these relations by starting from the fundamental solution to linear time-invariant (LTI) systems, then performing a discrete Fourier transform analytically, and finally projecting the result onto the SPOD modes at each frequency. Using the LTI solution as a starting point makes the method applicable to general linear systems, avoiding the requirement that the system be periodic, which might otherwise arise with a temporal Fourier basis. The operators involved are all linear and are precomputed, leaving only small matrix-vector multiplications to be done online. The method amounts to a space-time projection of a space-time solution operator, and we refer to it as \textit{Spectral Solution Operator Projection} (SSOP).
\\

We demonstrate the method on two problems: a linearized Ginzburg-Landau problem with a spatial dimension of $N_x = 220$ and an advection-diffusion problem with $N_x = 9604$. We show that, indeed, we can solve for the SPOD coefficients accurately, resulting in two-orders-of-magnitude lower error than even the projection of the solution onto the same number of POD modes, which is itself a lower bound for the error in any time-domain Petrov-Galerkin method, such as balanced truncation (BT) \citep{Moore81}. We show that this accuracy improvement does not come with an increase in CPU time; the computational cost of our method is similar to that of POD-Galerkin projection and balanced truncation, consistent with scaling estimates we derive.  
\\

The SPOD modes that form the basis for the method are obtained from data. A data-free version of the method is made possible by a connection between SPOD modes and the left singular vectors of the resolvent operator. The resolvent operator comes directly from the matrices that define the LTI system, so its singular value decomposition does not rely on data. As established by Ref. \citep{Towne2018spectral}, the left singular vectors of the resolvent operator at a given frequency are equivalent to the SPOD modes at the same frequency if the forcing in the system is spatially white. Many studies \citep{McKeon10,Beneddine16,Schmidt18,Nogueira19,Pickering20}, however, have shown that remarkable similarity between the two types of modes persists even when the forcing is far from white. Given that the SPOD modes and resolvent modes are similar, the latter serves as an excellent space-time basis for trajectories. We find that the data-free method that results from substituting the resolvent modes for the SPOD modes in the proposed method yields lower error than balanced truncation, the state-of-the-art data-free method.
\\

Though the space-time approach is uncommon, we are not the first to attempt it \citep{Yano14,Parish21,Choi21,Choi21b,Hoang22}. Previous methods have formed space-time basis vectors by assigning time dependence to POD modes, i.e., where each space-time basis vector is formed as the Kronocker product of a POD mode and a time dependent function. The time dependent functions are obtained from time series of the corresponding POD coefficients in the training data. We believe that the representational advantage of SPOD modes relative to previous choices of space-time basis, as well as their analytic time dependence, make them a more compelling choice for model reduction. 
\\

Using SPOD modes for linear model reduction has been explored before by Refs. \citep{Lin19} and \citep{Towne21}. Both methods can be viewed as reduced-order models for the frequency-domain equation $(i\omega \td{I} - \td{A}) \hat{\ket{q}}(\omega) = \td{B} \hat{\ket{f}}(\omega)$, where $\hat{\ket{q}}$ and $\hat{\ket{f}}$ are the Fourier transformed state and forcing, respectively, and $\td{A}$ and $\td{B}$ are the standard state-space matrices. However, this equation cannot capture transient behavior; it describes the steady state when the system is subjected to a periodic forcing, and this restriction precludes these ROMs from being applicable to many systems. In this paper, we derive a different frequency-domain equation that captures both transient and steady-state behavior and use it as the starting point for the proposed reduced-order model. By projecting a space-time solution operator onto SPOD modes, we prevent unretained modes from affecting the accuracy of retained coefficients. This was achieved by \citer{Towne21} (in the periodic case) by employing a Petrov-Galerkin projection, but not by \citer{Lin19} wherein a Galerkin projection was used. Another advancement of this work relative to \citers{Towne21, Lin19} is the use of different numbers of SPOD modes at different frequencies. This allows for more computational resources to be devoted to the more energetic frequencies. We also note that SPOD modes have been used to construct time-domain ROMs \citep{Chu21,Ji23}. However, these models are only superficially related to the space-time ROM proposed in this paper.
\\

Our approach may also be related to the harmonic balance method \citep{Hall02,Hall13}. In this technique, the governing equations for a temporally periodic system are Fourier-transformed in time, resulting in a set of nonlinear equations to be solved for the transformed fields. This has been applied to the periodic flows arising in turbomachinery \citep{Ekici07}, resulting in significant computational speedup due to the relatively small number of relevant harmonics. Harmonic balance does not employ a spatial reduction, and our method may be viewed as a spatially-reduced harmonic balance method for linear problems. Another important difference is that our method is applicable to non-periodic systems as well as periodic ones. 
\\

One requirement for the application of the proposed method bears mentioning here: the first step in the method is to take the temporal Fourier transform of the forcing. Therefore, the forcing on the entire time interval of interest must be available before starting computation in the reduced-order model. Any space-time or frequency-domain method that incorporates forcing, such as the ones mentioned above \citep{Lin19,Choi19,Parish21,Choi21,Choi21b,Yano14,Towne21}, is restricted in the same way. The proposed method is designed for applications where the forcing is known beforehand, such as, e.g., adjoint-based optimization and open loop control.
\\

The remainder of this paper is organized as follows. In Section~\ref{sec:POD_and_SPOD}, we discuss the properties of POD, space-time POD, and SPOD that are relevant to the method. We present the method in Section~\ref{sec:Methods}, and demonstrate it applied to a linearized Ginzburg-Landau problem and a scalar transport problem in Section~\ref{sec:results}. We conclude the paper in Section~\ref{sec:Conclusion}.

\begin{figure}
    \centering
    \begin{overpic}[scale = 0.65]{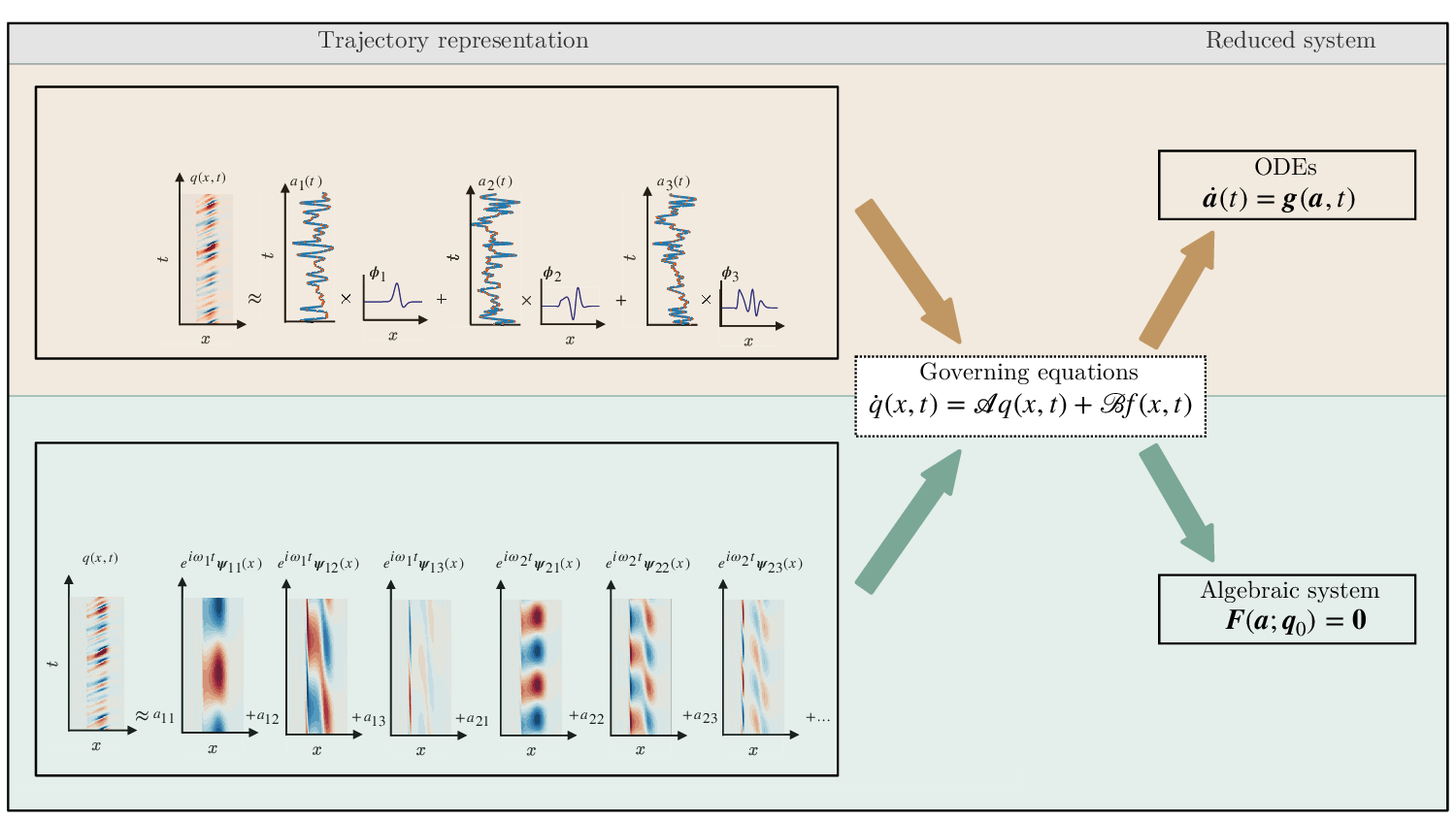}
    \put(17,47){$q(x,t) \approx \sum_{j = 1}^r a_j(t)\ket{\phi}_{j}(x)$}
    \put(13,22){$q(x,t) \approx \sum_{k = 0}^{N_\omega - 1} \sum_{j = 1}^{r_k} a_{kj} e^{i\omega_k t} \ket{\psi}_{kj}(x)$}

    \end{overpic}
    \caption{The proposed model reduction approach (bottom panel) compared against a standard space-only linear model reduction (top panel). To represent a trajectory, the space-only basis vectors are multiplied by time-dependent coefficients; the SPOD modes with their oscillatory time dependence are multiplied by static coefficients. In the space-only case, the coefficients are obtained by integrating a set of ODEs forward in time, whereas in the SPOD case, the coefficients are obtained by solving a linear algebraic system.}    \label{fig:graphical}
\end{figure}


\section{Space-only, space-time, and spectral POD} \label{sec:POD_and_SPOD}
We briefly review the space-only, space-time, and spectral forms of POD here; see Refs. \citep{Frame23} and \citep{Towne2018spectral} for additional details. The most significant point for the purposes of this paper is the fact that spectral POD modes approach space-time POD modes as the time interval becomes long, and thus are very efficient in representing trajectories.
\\

\subsection{Space-only POD}
Space-only POD aims to reconstruct snapshots of the state by adding together prominent modes weighted by expansion coefficients. In the continuous setting, the state $\vc{q}: \Omega \to \mathbb{C}^{N_v}$ is a function that maps elements of the spatial domain $\vd{x} \in \Omega$ to the $N_v$ state variables. The first space-only POD mode $\vc{\phi}_1: \Omega \to \mathbb{C}^{N_v}$ is defined to maximize the functional $\lambda[\vc{\phi}({\bf{x}})]$, which quantifies the expected value of the energy captured by its argument, 
\begin{subequations}
\begin{equation} \label{Sop:maximization}
    \lambda[\vc{\phi}({\bf{x}})] =  \frac{\mathbb{E}\big[ |\langle \vc{q}({\vd{x}}), \vc{\phi} ({\vd{x}})\rangle_{x} |^2 \big]}{\| \vc{\phi}({\vd{x}})\|_x^2 }  \text{,}
\end{equation}
\begin{equation} \label{Sop:maximization_with_functional}
    \vc{\phi}_1({\vd{x}})  = \argmax \lambda[\boldsymbol{\phi}({\bf{x}})] \text{.}
\end{equation}
\end{subequations}
The subsequent modes are defined to maximize the energy captured under the constraint that they are orthogonal to all previous ones,
\begin{equation} \label{Sop:maximization_with_functional_latter}
    \vc{\phi}_{j}(\vd{x}) = \argmax_{\langle \vc{\phi}({\vd{x}}), \vc{\phi}_{k<j} ({\vd{x}})\rangle_{x} = 0} \lambda[\boldsymbol{\phi}({\bf{x}})] \text{.}
\end{equation}
The space-only inner product $\langle \cdot, \cdot \rangle_{x}$, which defines the energy captured, is defined as an integral over the spatial domain $\Omega$,
\begin{equation} 
    \langle \vc{q}({\vd{x}}), \vc{\phi} ({\vd{x}})\rangle_{x} = \int_{\Omega} \vc{\phi}^* ({\vd{x}}) \tc{W}(\vd{x}) \vc{q}(\vd{x}) d\vd{x} \text{,}
\end{equation}
where $\tc{W}(\vd{x})$ is a weight matrix used to account for inter-variable importance or possibly to preference certain regions of the domain. One can show \citep{Towne2018spectral, Lumley67} that the solution to the optimization problem (\ref{Sop:maximization_with_functional},\ref{Sop:maximization_with_functional_latter}) is modes that are eigenfunctions of the space-only correlation tensor,
\begin{equation}
\int_{\Omega}\tc{C}(\vd{x}_1,\vd{x}_2) \tc{W}(\vd{x}_2) \vc{\phi}_j(\vd{x}_2) d\vd{x}_2 = \lambda_{j}\vc{\phi}_j(\vd{x}_1) \text{,}
\end{equation}
where the eigenvalue is equal to the energy of the mode, i.e., $\lambda_j = \lambda[\vc{\phi}_j]$, and the correlation tensor is 
\begin{equation}
    \tc{C}(\vd{x}_1,\vd{x}_2) = \mathbb{E}[\vc{q}(\vd{x}_1)\vc{q}^*(\vd{x}_2)] \text{.}
\end{equation}

\subsection{Space-time POD}
Whereas space-only POD modes optimally represent snapshots, space-time POD modes optimally represent trajectories over the time window $[0,T]$. The formulation is much the same as in space-only POD; the modes optimize the expected energy (\ref{Sop:maximization_with_functional},\ref{Sop:maximization_with_functional_latter}), but in space-time POD, the inner product involved in defining the energy functional includes time as well as space,
\begin{equation} \label{eq:ip_pod}
    \langle \vc{q}({\vd{x}},t), \vc{\phi} ({\vd{x}},t)\rangle_{x,t} = \int_{0}^{T} \int_{\Omega} \vc{\phi}^* ({\vd{x}},t) \tc{W}(\vd{x}) \vc{q}(\vd{x},t) d\vd{x} \ dt \text{.}
\end{equation}
The space-time POD modes that solve this optimization are eigenfunctions of the space-time correlation $\tc{C}(\vd{x}_1,t_1,\vd{x}_2,t_2) =  \mathbb{E}[\vc{q}(\vd{x}_1,t_1)\vc{q}^*(\vd{x}_2,t_2)]$, and the eigenvalues represent the energy of each mode. The most important property of space-time POD modes for the purpose of this paper is that the space-time POD reconstruction of a trajectory achieves lower error, on average, than the reconstruction with the same number of modes in any other space-time basis. More concretely, using the first $r$ space-time POD modes to reconstruct the trajectory,
\begin{equation}
    \tilde{\vc{q}}(\vd{x},t) = \sum_{j = 1}^{r} \vc{\phi}_j ({\vd{x}},t) \langle \vc{q}(\vd{x},t) , \vc{\phi}_j ({\vd{x}},t) \rangle_{x,t} \text{,}
\end{equation}
yields lower expected error
\begin{equation} \label{eq:exerr_pod}
    \mathbb{E} [\|\tilde{\vc{q}}(\vd{x},t) - \vc{q}(\vd{x},t)\|_{x,t}^2 ] 
\end{equation}
than would any other space-time basis. The expected error \eqref{eq:exerr_pod} is measured over space and time using the norm $\| \cdot \|_{x,t}$ induced by the inner product \eqref{eq:ip_pod}.
\\

\subsection{Spectral POD}
Spectral POD is most easily understood as the frequency domain variant of space-only POD for statistically stationary systems. In other words, SPOD modes at a particular frequency optimally reconstruct (in the same sense as above) the state at that frequency, on average. The property that makes them attractive for model reduction, however, is that SPOD modes are also the long-time limit of space-time POD modes for statistically stationary systems. These ideas are made precise below, but for a more complete discussion, see Ref. \citep{Towne2018spectral}.
\\

Spectral POD modes at frequency $k$ maximize
\begin{equation} \label{eq:POD:spectral_op}
     \lambda_k[\vc{\psi}({\bf{x}})] =  \frac{\mathbb{E}\big[ |\langle \hat{\vc{q}}_k({\vd{x}}), \vc{\psi} ({\vd{x}})\rangle_{x} |^2 \big]}{\| \vc{\psi}({\vd{x}})\|_x^2 }  \text{,}
\end{equation}
again subject to the constraint that each mode $\vc{\psi}_{k,j}(\vd{x})$ is orthogonal to the previous ones at that frequency $\vc{\psi}_{k,i<j}(\vd{x})$. The Fourier-transformed state $\hat{\vc q}_k: \Omega \to \mathbb{C}^{N_v}$ is defined as 
\begin{equation}
    \hat{\vc q}_k(\vd{x}) = \int_{-\infty}^{\infty} e^{-i\omega_k t}\vc{q}(\vd{x},t) dt \text{.}
\end{equation}
The solution to the optimization (\ref{eq:POD:spectral_op}) is that the modes are eigenvectors of the cross-spectral density
   $\tc{S}_{k}(\vd{x}_1,\vd{x}_2) = \mathbb{E}[\hat{\vc q}_k(\vd{x}_1) \hat{\vc q}_k^*(\vd{x}_2)]$,
\begin{equation}
    \int \limits_{\Omega}\tc{S}_{k}({\vd{x}}_1,{\vc{x}}_2)\tc{W}({\vd{x}}_2)\vc{\psi}_{k,j}({\vd{x}}_2)d{\vd{x}}_2 = \lambda_{k,j}\vc{\psi}_{k,j}({\vd{x}}_1) \text{.}
\end{equation}
The eigenvalue is again equal to the energy of the mode, i.e., $\lambda_{k,j} = \lambda_{k}[\vc{\psi}_{k,j}]$. Modes at frequency $k$ have an implicit time dependence of $e^{i\omega_k t}$.
\\

SPOD modes and their energies become identical to space-time POD modes as the time interval on which the latter are defined becomes long \citep{Lumley67,Towne2018spectral,Frame23}. Thus, the SPOD modes with the largest energies among all frequencies are the dominant space-time POD modes (for long times) and are most efficient for reconstructing long-time trajectories. We denote by $\tilde{\lambda}_j$ the $j$-th largest SPOD eigenvalue among all frequencies, which may be compared to the space-time POD eigenvalues. The convergence in the energy of the trajectory they capture is relatively fast, so for time intervals beyond a few correlation times, the SPOD modes capture nearly as much energy as the space-time modes \citep{Frame23}. If the simulation time of a reduced-order model is long enough for this convergence to be met, the ability of the SPOD modes to capture structures is not diminished relative to that of space-time POD modes. 
\\

SPOD modes also have two properties that make them more suitable for model reduction than space-time modes: they have analytic time dependence, and they are separable in space and time. The former makes some analytic progress possible in writing the equations that govern the modes and enables Fourier theory to be applied. The latter means that storing the modes requires $N_t$ times less memory, where $N_t$ is the number of times steps in the simulation.
\\

\begin{figure}
    \centering
    \input{Figures/DOF_vs_T}
    \includegraphics{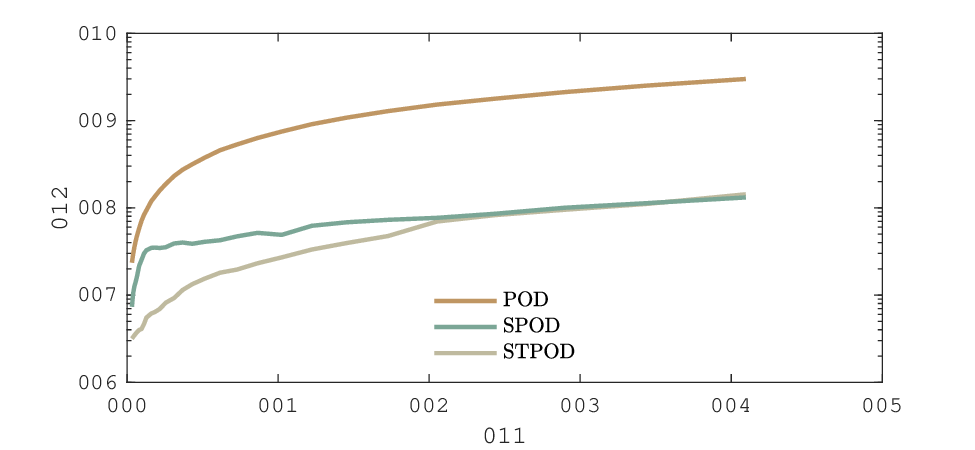}
    \caption{Number of degrees of freedom (DOFs) required to achieve $98\%$ representation accuracy of trajectories as a function of the length of the time interval of the trajectory $[0,T]$. For POD, one must specify all the mode coefficients at every time step, whereas spectral and space-time POD modes are themselves time-dependent. Thus, by leveraging spatiotemporal correlations, fewer DOFs are needed to represent a trajectory to a given accuracy by specifying the SPOD or space-time POD coefficients. As the time interval becomes long, the SPOD and space-time POD modes become equally efficient at representing trajectories. The data come from the Ginzburg-Landau system introduced in Section~\ref{sec:results}, and the details on the comparison of POD and SPOD with the same number of degrees of freedom can be found there.}
    \label{fig:DOF_vs_T}
\end{figure}

Figure~\ref{fig:DOF_vs_T} shows the convergence in the representational ability of space-time POD and SPOD as the time interval becomes long. Specifically, to represent trajectories with some level of accuracy, $98 \%$, say, one needs the same number of SPOD coefficients as space-time POD coefficients if the interval is long compared to the correlation time in the system. This convergence in representation ability occurs because the SPOD modes themselves converge to space-time POD modes in the limit of a long time interval. With space-only POD, one must specify the coefficients for every time step, which leads to a far less efficient encoding of the data because the coefficients are highly correlated from one time step to the next. That SPOD modes are near-optimal in representing trajectories, and that they are substantially more efficient than space-only POD modes, motivate this work. If one can efficiently solve for some number of the SPOD coefficients of a trajectory, then these coefficients will lead to substantially lower error than solving for the same number of space-only POD coefficients.

\subsection{Discretization and truncation}
Upon spatial discretization, the state and modes become vectors in $\mathbb{C}^{N_x}$, where $N_x$ denotes the number of gridpoints multiplied by $N_v$. Frequency is also discretized, and a finite number $N_\omega$ of evenly spaced frequencies is retained. The lowest one, $\omega_1$, induces a time $T = 2\pi/\omega_1$, which determines the interval $[0,T]$ on which the modes are periodic. The trajectories themselves are, of course, not periodic on this interval, so $T$ is the longest we may use SPOD modes for prediction, though, if a longer prediction is needed, the method may be repeated. The fastest frequency corresponds to a time step $\Delta t$. The discrete Fourier transform of a trajectory $\ket{q}(t)$ at frequency $\omega_k = 2\pi k /T$ is defined as 
\begin{equation} \label{eq:DFT}
    \hat{\ket{q}}_k = \sum_{j = 0}^{N_\omega - 1} \ket{q}(j \Delta t) e^{-i\omega_k j \Delta t} \text{,}
\end{equation}
where $\hat{\ket{q}}_k $ and $\ket{q}(j \Delta t)$ are both in $\mathbb{C}^{N_x}$.
\\

With the spatial discretization, integration over space becomes matrix-vector multiplication. For example, the full set of discrete SPOD modes at frequency $\omega_k$ are defined as the eigenvectors of the (weighted) cross-spectral density matrix 
\begin{equation}
    \td{S}_k \td{W} \td{\Psi}_{k,\text{full}} = \td{\Psi}_{k,\text{full}}\td{\Lambda}_{k,\text{full}} \text{.}
\end{equation}
Here, $\td{S}_k \in \mathbb{C}^{N_x \times N_x}$ is the cross-spectral density, and $\td{W} \in \mathbb{C}^{N_x \times N_x}$ is a weight matrix. $\td{\Psi}_{k,\text{full}} \in \mathbb{C}^{N_x \times N_x}$ is the matrix with the full set discrete SPOD modes as its columns, and $\td{\Lambda}_{k,\text{full}}$ is the diagonal matrix of corresponding eigenvalues, which are the SPOD mode energies. However, in practice, the matrix $\td{S}_k$ is not formed, and too little data will be available to create a full-rank set of SPOD modes. Instead, the SPOD modes are calculated by the method of snapshots \citep{Sirovich87,Towne2018spectral} or by using the singular value decomposition. In particular, given $r_d$ trajectories from which to obtain SPOD modes, each $N_{\omega}$ time steps in length, the discrete Fourier transform (DFT) of each trajectory is taken. This yields $r_d$ realizations of the $k$-th frequency, for every frequency $k$. These can be formed into a data matrix
\begin{equation}
    \td{Q}_k = [\hat{\ket{q}}_k^1, \hat{\ket{q}}_k^2, \dots,  \hat{\ket{q}}_k^{r_d}] \text{,}
\end{equation}
where $\hat{\ket{q}}_k^i \in \mathbb{C}^{N_x}$ is the $k$-th frequency of the DFT of the $i$-th trajectory. The SPOD modes at frequency $\omega_k$ and the associated energies may then be obtained by first taking the singular value decomposition $\td{U}\td{\Sigma}\td{V}^* = 1/\sqrt{r_d}\td{W}^{1/2}\td{Q}_k$. The $r_d$ available SPOD modes $\td{\Psi}_k^{r_d} \in \mathbb{C}^{N_x \times r_d}$ are then given by $\td{W}^{-1/2}\td{U}$ and the energies $\td{\Lambda}_k^{r_d} \in \mathbb{R}^{r_d \times r_d}$ by $\td{\Sigma}^2$ \citep{Towne2018spectral}. In practice, the $r_d$ trajectories are usually obtained as (possibly overlapping) sub-trajectories of one long trajectory, as is done in Welch's method \citep{Welch67} for computing power spectra.
\\

The set of available SPOD modes at the $k$-th frequency $\td{\Psi}_k^{r_d}$ must be truncated in forming a reduced-order model, however, the order of this truncation should not be the same for all frequencies. Intuitively, allocating more SPOD modes to the energetic frequencies leads to more accuracy than keeping the number of modes constant across frequency. We denote the mean number of modes retained at each frequency as $r$, thus $N_\omega r$ modes are retained in total. We determine the number of SPOD modes retained at frequency $\omega_k$ as the number of modes at this frequency that are among the $N_\omega r$ most energetic overall. That is,
\begin{equation}
    r_k = | \{ l : \lambda_{k,l} \geq \tilde{\lambda}_{N_\omega r} \} | \text{,}
\end{equation}
where  $\tilde{\lambda}_{i}$ denotes the $i$-th largest eigenvalue over all frequencies. With this notation established, we denote the retained SPOD modes at frequency $\omega_k$ as $\td{\Psi}_k \in \mathbb{C}^{N_x \times r_k}$, and the corresponding energies as $\td{\Lambda}_k \in \mathbb{R}^{r_k \times r_k}$
\\

The trajectory $\ket{q}(t)$ on the interval $[0,T]$ is then approximated using the retained SPOD modes as
\begin{equation}
    \ket{q}(t) \approx \frac{1}{N_\omega}\sum_{k = 0}^{N_{\omega} - 1} \td{\Psi}_k \ket{a}_k e^{i\omega_k t} \text{.}
\end{equation}
The vector of expansion coefficients at frequency $\omega_k$ is given by 
\begin{equation} \label{eq:Uncor:a from qhat}
    \ket{a}_k = \td{\Psi}_k^* \td{W} \hat{\ket{q}}_k \in \mathbb{C}^{r_k} \text{.}
\end{equation}

\section{Spectral solution operator projection method} \label{sec:Methods}
Our goal is to derive an SPOD-based method to solve the linear ordinary differential equation
\begin{subequations}\label{eq:ODE}
\begin{equation} 
    \dot{\ket q}(t) = \td{A} \ket{q}(t) + \td{B}\ket{f}(t) \text{,}
\end{equation}
\begin{equation}
    \ket{y}(t) = \td{C} \ket{q}(t)
\end{equation}
\end{subequations}
on the interval $t \in [0,T]$, where $\ket{q}(t) \in \mathbb{C}^{N_x}$ is the state, $\td{A} \in \mathbb{C}^{N_x \times N_x}$ is the system matrix, $\ket{f}(t) \in \mathbb{C}^{N_f}$ is some known forcing that is mapped onto the system by the matrix $\td{B} \in \mathbb{C}^{N_x \times N_f}$, and $\ket{y} \in \mathbb{C}^{N_y}$ is an observable extracted from the state by the matrix $C \in \mathbb{C}^{N_y \times N_x}$. Given the forcing and initial condition $\ket{q}(0) = \ket{q}_0$, our goal is to find the retained SPOD coefficients for the trajectory $\ket{q}(t)$, thereby obtaining the near-optimal rank-$N_\omega r$ space-time representation of the trajectory. With these coefficients, $\ket{y}(t)$ can be easily obtained taking the inverse DFT of $\hat{\ket{y}}_k = \td{C}\td{\Psi}_k \ket{a}_k$. 
\\

The starting point for finding the retained coefficients is \eqref{eq:Uncor:a from qhat}, which gives $\ket{a}_k$ in terms of $\hat{ \ket{q}}_k$. The Fourier-transformed state $\hat{\ket q}_k$ must be obtained from the known forcing and initial condition, and we derive this relation in the following subsection.

\subsection{Frequency-domain equation}
We begin by inserting the analytic solution to \eqref{eq:ODE} into the definition of the DFT \eqref{eq:DFT},
\begin{equation} \label{eq:DFT_of_TD}
      \vdh{q}_k = \sum_{j = 0}^{N_\omega - 1} e^{-i\omega_k j \Delta t}\left( e^{\td{A}j \Delta t}\ket{q}_0 + \int_0^{j \Delta t} e^{\td{A} (j \Delta t - t')} \td{B} \ket{f}(t') \ dt' \right) \text{.}
\end{equation}
The term in parentheses is the well-known solution to the linear ODE \eqref{eq:ODE} \citep{Grigoriu21}. Analytic progress can be made with the assumptions that the forcing can be written as $\ket{f}(t) = \sum_{k = 0}^{N_\omega - 1} \vdh{f}_k e^{i\omega_k t}$, and that the matrices $i \omega_k \td{I} - \td{A}$ are invertible for all $k$. To aid in the ensuing discussion, we split \eqref{eq:DFT_of_TD} into two parts: the response to the initial condition and the response to the forcing. The first of these is 
\begin{equation} \label{just IC term DFT}
    \hat{\ket{q}}_{k,ic} = \sum_{j = 0}^{N_\omega - 1}  e^{(\td{A} - i \omega_k \td{I})j\Delta t} \ket{q}_0 \text{.}
\end{equation}
This term may be evaluated by noticing that it is a matrix geometric sum, so, it can be written as 
\begin{equation}
   \hat{\ket{q}}_{k,ic} = (\td{I} - e^{(\td{A} - i \omega_k \td{I})\Delta t})^{-1} (\td{I} - e^{(\td{A} - i \omega_k)N_\omega\Delta t}) \ket{q}_0 \text{.}
\end{equation}
Because $e^{i\omega_k N_\omega \Delta t} = 1$, this simplifies to
\begin{equation}
    \hat{\ket{q}}_{k,ic} = (\td{I} - e^{(\td{A} - i \omega_k)\Delta t})^{-1} (\td{I} - e^{\td{A}T})\ket{q}_0 \text{.}
\end{equation}
Note that the assumption that $i \omega_k \td{I} - \td{A}$ is invertible implies that $\td{I} - e^{(\td{A} - i \omega_k \td{I})\Delta t}$ is also invertible. 
\\

Second, the response to the forcing is
\begin{equation}
    \sum_{j = 0}^{N_\omega - 1} e^{-i\omega_k j \Delta t} \int_0^{j \Delta t} e^{\td{A} (j \Delta t - t')} \td{B} \ket{f}(t') \ dt' \text{.}
\end{equation}
To evaluate this term, we first insert the Fourier expansion of the forcing and take the constant matrix exponential out of the integral, giving
\begin{equation}
   \hat{\ket{q}}_{k,force} = \frac{1}{N_\omega} \sum_{j = 0}^{N_\omega - 1} e^{-i\omega_k j \Delta t} e^{\td{A} j \Delta t} \int_0^{j \Delta t} \sum_{l = 0}^{N_\omega - 1}  e^{(i\omega_l\td{I} - \td{A})t'} \td{B} \hat{\ket{f}}_l \ dt' \text{.}
\end{equation}
Integration inverts the matrix in the exponential, and the expression evaluates to
\begin{equation} \label{eq:intermediate3}
    \hat{\ket{q}}_{k,force} = \frac{1}{N_\omega} \sum_{j = 0}^{N_\omega - 1} e^{-i\omega_k j \Delta t} \sum_{l = 0}^{N_\omega - 1} \td{R}_l \left( e^{i\omega_l j \Delta t} - e^{\td{A} j \Delta t} \right) \td{B} \hat{\ket{f}}_l \text{.}
\end{equation}
We refer to $\td{R}_k = (i\omega_k \td{I} - \td{A})^{-1}$ as the resolvent operator. Equation \eqref{eq:intermediate3} can then be written as
\begin{equation} \label{eq:intermediate4}
    \hat{\ket{q}}_{k,force} = \frac{1}{N_\omega} \sum_{j = 0}^{N_\omega - 1}  \sum_{l = 0}^{N_\omega - 1} \td{R}_l \left( e^{i(\omega_l - \omega_k) j \Delta t} - e^{(\td{A} - i\omega_k \td{I}) j \Delta t} \right) \td{B} \hat{\ket{f}}_l \text{.}
\end{equation}
The frequency difference term in \eqref{eq:intermediate4} evaluates to zero for $\omega_l \neq \omega_k$, and the other term may be evaluated by the same geometric sum argument as before. The entire forcing term in \eqref{eq:DFT_of_TD} becomes
\begin{equation}
    \hat{\ket{q}}_{k,force} = \td{R}_k \td{B} \hat{\ket{f}}_k - \frac{1}{N_\omega} (\td{I} - e^{(\td{A} - i \omega_k)\Delta t})^{-1} (\td{I} - e^{\td{A}T}) \sum_{l = 0}^{N_\omega - 1} \td{R}_l \td{B} \hat{\ket{f}}_l \text{.}
\end{equation}
Recombining our simplified expressions for $\hat{\ket{q}}_{k,ic}$ and $\hat{\ket{q}}_{k,force}$, we obtain the following equation for the $k$-th component of the DFT of the state:
\begin{equation} \label{solution}
    \hat{\ket{q}}_k = \td{R}_k \td{B}  \hat{\ket{f}}_k + \left(\td{I} -e^{(\td{A} - i \omega_k)\Delta t} \right)^{-1} \left(\td{I} - e^{\td{A}T} \right) \left(\ket{q}_0 - \frac{1}{N_\omega} \sum_{l = 0}^{N_\omega - 1} \td{R}_l \td{B} \hat{\ket{f}}_l \right) \text{.}
\end{equation}
\\

The first term in \eqref{solution} $\td{R}_k \td{B}  \hat{\ket{f}}_k$, is familiar: it is the component at frequency $\omega_k$ of the steady-state response to a periodic forcing. The second term in \eqref{solution} represents the transient. This transient may persist for the entirety of the interval $[0,T]$, so including it is crucial. 
\\

\citers{Lin19, Towne21} took an expression equivalent to the first term as their starting point for SPOD-based model reduction. Specifically, they began with the equation 
\begin{equation} \label{eq:FD_wrong}
    \td{L}_k \vdh{q}_k = \td{B} \vdh{f}_k \text{,}
\end{equation}
where $\td{L}_k = (i\omega_k \td{I} - \td{A})$ is the inverse of the resolvent operator. We stress that unless the solution $\ket{q}(t)$ is $T$-periodic, this relation is incorrect, and this issue was not fully appreciated in either \citer{Lin19} or \citer{Towne21}. We note that \citep{Martini19} has discussed this issue in another context.
\\

\citer{Lin19} proceeded to reduce \ref{eq:FD_wrong} by applying a Galerkin projection with $\td{\Psi}_k$ as the trial basis, resulting in the equation $\ket{a}_k = (\td{\Psi}_k^{*}\td{W}\td{L}_k\td{\Psi}_k^{})^{-1} \td{B} \vdh{f}_k$. Even in $T$-periodic systems, this reduction does not recover the exact SPOD coefficients --- the unretained modes influence the retained coefficients. An innovation made \citer{Towne21} was to use a Petrov-Galerkin projection of \eqref{eq:FD_wrong}, resulting in the equation $\ket{a}_k = (\td{\Phi}_k^{F*}\td{W}\td{L}_k\td{\Psi}_k^{})^{-1}\td{\Phi}_k^{F*}\td{W}\td{B}\vdh{f}$. \citer{Towne21} showed that the test basis $\td{\Phi}_k^{F}$ can be chosen in such a way that the retained SPOD coefficients are exact (for $T$-periodic systems). The expression for the test basis $\td{\Phi}_k^{F}$ required to achieve this involved the statistics of the forcing and was written in terms of the singular value decomposition of the product of the resolvent operator and a matrix that contained these forcing statistics.
\\

The $T$-periodic Petrov-Galerkin method from \citer{Towne21} may be alternatively derived by left-multiplying the full-order equation $\vdh{q}_k = \td{R}_k\td{B} \vdh{f}_k$ by $\td{\Psi}^*_k \td{W}$, resulting in the following relation for the SPOD coefficients $\ket{a}_k = \td{\Psi}^*_k \td{W}\td{R}_k\td{B} \vdh{f}_k$. Since our goal is to arrive at a ROM without the constraint of periodicity, we instead apply the left-multiplication to \eqref{solution}.
\\

While the frequency domain equation \eqref{solution} remains valid regardless of the stability of the system as long as the assumptions listed above are met, we do not recommend the method in the case of an unstable system unless $T\lambda_{max}$ is relatively close to unity. When there is significant exponential growth along the unstable eigenvectors, i.e., $e^{T\lambda_{max}}$ is large, capturing these directions, and the projection of the forcing onto them is crucial, so eigensystem methods are superior. 
\subsection{Reduction and operator approximations}
As described above, we proceed by left-multiplying \eqref{solution} according to \eqref{eq:Uncor:a from qhat},
\begin{equation} \label{eq:corr_pgp_toapprox}
        \ket{a}_k = \td{\Psi}^{*}_k \td{W} \td{R}_k \td{B} \hat{\ket{f}}_k +  \td{\Psi}^{*}_k \td{W} \left(\td{I} -e^{(\td{A} - i \omega_k \td{I})\Delta t} \right)^{-1} \left(\td{I} - e^{\td{A}T} \right) \left(\ket{q}_0 - \frac{1}{N_\omega} \sum_{l = 0}^{N_\omega - 1} \td{R}_l \td{B} \hat{\ket{f}}_l \right) \text{.}
\end{equation}
If the system is small, all operators may be computed analytically. However, for large systems, direct computation of, e.g., the inverse that defines the resolvent operator or the matrix exponentials, is not tractable. For these systems, these operators must be approximated.
\\

\subsubsection{Steady-state operator}
We begin with the operator $\td{M}_k = \td{\Psi}^{*}_k \td{W} \td{R}_k \td{B}$. An accurate and simple-to-implement approximation of $\td{M}_k$ can be obtained by leveraging the availability of data as follows. Defining
\begin{equation} \label{eq:gdef}
    \ket{g}_k^i  = \td{L}_k\hat{\ket q}_k^i \text{,}
\end{equation}
where, again, $\td{L} =  (i\omega_k \td{I} - \td{A})$ is the inverse of the resolvent and $\hat{\ket q}_k^i$ is the $i$-th realization of $\hat{\ket q}_k$ in the training data, we have
\begin{equation}
    \hat{\ket q}_k^i = \td{R}_k \ket{g}_k^i \text{.}
\end{equation}
Forming each $\ket{g}_k^i$ is not computationally expensive so long as $\td{A}$ is sparse. Using the many realizations of the training data (the same ones used to generate the SPOD modes), we have
\begin{equation}
    \td{Q}_k = \td{R}_k \td{G}_k \text{,}
\end{equation}
where $\td{Q}_k = [\hat{\ket{q}}_k^1, \hat{\ket{q}}_k^2, \dots, \hat{\ket{q}}_k^{r_d}] $ and $\td{G}_k = [\ket{g}_k^1, \ket{g}_k^2, \dots, \ket{g}_k^{r_d}]$. We approximate the resolvent operator $\td{R}_k$ as
\begin{equation}
    \td{R}_k \approx \td{Q}_k \td{G}^+ \text{,}
\end{equation}
where the pseudoinverse is defined $\td{G}^+ = (\td{G}^* \td{W} \td{G})^{-1} \td{G}^* \td{W}$. It may be shown that the action of this approximate resolvent is equivalent to an orthogonal projection into the column space of $\td{G}_k$, followed by a multiplication by the resolvent. In other words, $\td{Q}_k \td{G}^+ = \td{R}_k \td{P}_g$, where $\td{P}_g$ is an orthogonal projection matrix (in the $\td{W}$-based norm). This means that $\td{Q}_k \td{G}^+$ is an effective approximation of $\td{R}_k$ to the extent that the vectors to which it is applied are near the column space of $\td{G}_k$. By inserting this approximation of the resolvent operator into the expression $\td{M}_k = \td{\Psi}^{*}_k \td{W} \td{R}_k \td{B}$, we define the matrix $\td{E}_k$ as
\begin{equation}
    \td{M}_k \approx \td{E}_k = \td{\Psi}^{*}_k \td{W} \td{Q}_k \td{G}^+_k \td{B}  \text{.}
\end{equation}
The total cost of these operations scales linearly with $N_x$ and quadratically with $r_d$, thus the approximation avoids the superlinear $N_x$ scaling required by most methods for computing the (non-approximate) action of the resolvent operator. We note that if more accuracy is required, time-stepping approaches, such as those developed by \cite{Farghadan25}, may be used.  

\subsubsection{Transient operator} \label{sec:transient_op}
Next, we approximate the operator $ \td{\Psi}^{*}_k \td{W} \left(\td{I} -e^{(\td{A} - i \omega_k \td{I})\Delta t} \right)^{-1} \left(\td{I} - e^{\td{A}T} \right)$. This is accomplished using the $r_d$ available SPOD modes at each frequency. We first multiply by the identity, expressed as $\td{\Psi}_{k,\text{full}} \td{\Psi}^{*}_{k,\text{full}} \td{W}$, in various places,
\begin{subequations}
\begin{align}
\begin{split} \label{eq:corr_pgp:intermediate1}
\td{\Psi}^{*}_k \td{W} & \left(\td{I} -e^{(\td{A} - i \omega_k \td{I})\Delta t} \right)^{-1} \left(\td{I} - e^{\td{A}T} \right)\\
&= \td{\Psi}^{*}_k \td{W} \left(\td{I} -e^{(\td{A} - i \omega_k)\Delta t} \right)^{-1} \td{\Psi}_{k,\text{full}} \td{\Psi}^{*}_{k,\text{full}} \td{W} \left(\td{I} - e^{\td{A}T} \right) \td{\Psi}_{k,\text{full}} \td{\Psi}^{*}_{k,\text{full}} \td{W} \text{.}
\end{split}\\
\begin{split}  \label{eq:corr_pgp:intermediate2}
&=\td{P}_k \left(\td{I} -e^{(\td{\Psi}^{*}_{k,\text{full}} \td{W} \td{A}\td{\Psi}_{k,\text{full}} - i \omega_k \td{I})\Delta t} \right)^{-1}  \left(\td{I} - e^{\td{\Psi}^{*}_{k,\text{full}} \td{W} \td{A}\td{\Psi}_{k,\text{full}}T} \right)  \td{\Psi}^{*}_{k,\text{full}} \td{W} \text{.}
\end{split}
\end{align}
\end{subequations}
The full-rank set of SPOD modes are brought into the matrix exponentials in \eqref{eq:corr_pgp:intermediate2}, and it is straightforward to show that this does not introduce any approximation. In \eqref{eq:corr_pgp:intermediate2}, the matrix $\td{P}_k = \begin{bmatrix}
    \td{I}_{r_k} && \td{0}
\end{bmatrix} \in \mathbb{R}^{r_k \times r_d}$ selects the first $r_k$ rows of the matrix it multiplies. Finally, truncating the operators in \eqref{eq:corr_pgp:intermediate2}, i.e., $\td{\Psi}_{k,\text{full}} \to \td{\Psi}_k^{r_d}$, and denoting $\tilde{\td{A}} = \td{\Psi}^{r_d*}_k \td{W} \td{A}\td{\Psi}^{r_d}_k$, the approximated term is 
\begin{equation}
    \td{\Psi}^{*}_k \td{W} \left(\td{I} -e^{(\td{A} - i \omega_k \td{I})\Delta t} \right)^{-1} \left(\td{I} - e^{\td{A}T} \right) \approx \td{P}_k \left(\td{I} -e^{(\tilde{\td{A}}_k - i \omega_k \td{I})\Delta t} \right)^{-1} \left(\td{I} - e^{\tilde{\td{A}}_kT} \right) \td{\Psi}^{r_d*}_k \td{W} \text{.}
\end{equation}
As desired, all matrix exponentials and inverses are of size $r_d \times r_d$, which makes them tractable.
\\

The forcing sum in \eqref{eq:corr_pgp_toapprox}, which is difficult to compute directly because it involves resolvents, must also be approximated. In the unreduced equations, this term is the same at each frequency, so the sum over frequencies only needs to be computed once. Any approximation of this term should be the same for each frequency to avoid quadratic scaling in $N_\omega$. The natural choice is to approximate each term by 
\begin{equation}
    \td{R}_l \td{B} \hat{\ket{f}}_l \approx \td{\Psi}_l^{}\td{\Psi}_l^{*} \td{W} \td{R}_l \td{B} \hat{\ket{f}}_l \text{.}
\end{equation}
This approximation is accurate because the SPOD modes at the $l$-th frequency are the best basis for $\hat{\ket{q}}_l$ and are thus a very good basis for $\td{R}_l \td{B} \hat{\ket{f}}_l$, which is the steady-state component of $\hat{\ket{q}}_l$. The operator $\td{\Psi}_l^{*} \td{W} \td{R}_l \td{B}$ may again be approximated with $\td{E}_l$. 
\\

The equations can, at this point, be written as 
\begin{equation} \label{eq:simpler_pgp_noib}
     \ket{a}_k = \td{E}_k \hat{\ket{f}}_k + \td{F}_k  \left(\ket{q}_0 - \frac{1}{N_\omega} \sum_{l} \td{\Psi}_l \td{E}_l \hat{\ket{f}}_l \right) \text{,}
\end{equation}
where $\td{F}_k = \left(\td{I} -e^{(\tilde{\td{A}}_k - i \omega_k \td{I})\Delta t} \right)^{-1} \left(\td{I} - e^{\tilde{\td{A}}_kT} \right) \td{\Psi}^{*}_k \td{W} \in \mathbb{C}^{r_k \times N_x}$. The operators have been approximated, so far, to avoid  $\mathcal{O}(N_x^3)$ scaling in the offline phase of the algorithm. To avoid $\mathcal{O}(N_x)$ scaling in the online phase, one final approximation must be made. The term in parentheses in (\ref{eq:simpler_pgp_noib}) is multiplied on the left by $\td{F}_k$ for every frequency $\omega_k$, leading to $N_x N_\omega r$ scaling. This can be avoided by storing the term in parentheses in (\ref{eq:simpler_pgp_noib}) in a rank-$p$ reduced basis $\td{\Phi} \in \mathbb{C}^{N_x \times p}$ and precomputing the product of this basis with each $\td{F}_k$. This basis should represent the initial condition and forcing sum terms accurately, and in practice, we choose POD modes of the state. With this approximation, the equations become
\begin{equation} \label{eq:corr_pgp:final_eq}
     \ket{a}_k = \td{E}_k \hat{\ket{f}}_k + \td{H}_k   \left(\td{\Phi}^*\td{W} \ket{q}_0 - \frac{1}{N_\omega} \sum_{l} \td{T}_l \td{E}_l \hat{\ket{f}}_l \right) \text{,}
\end{equation}
where $\td{H}_k = \td{F}_k \td{\Phi} \in \mathbb{C}^{r_k \times p}$ and $\td{T}_l = \td{\Phi}^* \td{W} \td{\Psi}_l \in \mathbb{C}^{p \times r_l}$. Given an initial condition $\ket{q}_0$ and forcing $\ket{f}(t)$, the online stage of the method consists of taking the Fourier transform of the forcing, inserting this and the initial condition into \eqref{eq:corr_pgp:final_eq} in order to get the SPOD coefficients, then transforming back to the time domain. The details are given in the following subsection.

\subsection{Formal statement of the algorithm and scaling}
The offline and online phases of the SSOP method are shown in Algorithms~\ref{alg:offline} and \ref{alg:online}, respectively. 
\begin{algorithm} 
\caption{SSOP (offline)}\label{alg:offline}
\begin{algorithmic}[1]
\State \textbf{Inputs:} $\td{A}$, $\td{B}$, $\td{C}$, $\td{W}$, $\{ \td{\Psi}^{r_d}_i \}$, $\{ \td{\Lambda}^{r_d}_i \}$,$\{ \td{G}_i\}$, $\{ r_i \}$, $\td{\Phi}$
\For{ $k \in \{1,2,\dots, N_{\omega} \}$}
\State $\td{\Psi}_k^{} \gets [\ket{\psi}_{k,1}, \dots , \ket{\psi}_{k,r_k}]$ \Comment{Retained SPOD modes}
\State $\td{L}_k \gets i\omega_k \td{I} - \td{A}$ \Comment{Inverse of resolvent}

\State $\td{E}_k \gets  \td{\Psi}^{*}_k \td{W} \td{Q}_k \td{G}^+_k \td{B}$ \Comment{Precomputation of first operator}

\State $\tilde{\td{A}}_k \gets \td{\Psi}^{r_d*}_k \td{W} \td{A} \td{\Psi}_k^{r_d}$ \Comment{Reduced $\td{A}$ to compute matrix exponentials}
\State $\td{P}_k \gets \begin{bmatrix}
    \td{I}_{r_k} && \td{0}
\end{bmatrix}$ \Comment{Row selector matrix}
\State $\td{H}_k \gets \td{P}_k (\td{I} - e^{(\tilde{\td{A}}_k - i\omega_k \td{I})\Delta t} )^{-1} (\td{I} - e^{\tilde{\td{A}}_k T}) (\td{\Psi}_k^{r_d*} \td{W} \td{\Phi}$) \Comment{Precomputation of second operator}
\State $\td{T}_k \gets \td{\Phi}^* \td{W} \td{\Psi}_k$ \Comment{Precomputation of third operator}
\State $\td{C}^{\td{\Psi}}_k \gets \td{C}\td{\Psi}_k$ \Comment{$\td{C}$ in SPOD basis}

\EndFor

\end{algorithmic}
\vspace{1\baselineskip}
\textbf{Inputs:}  $\td{A}$, $\td{B}$, $\td{C}$, the system matrices; $\td{W}$, the weight matrix; $\{ \td{\Psi}^{r_d}_i \}$, the $r_d$ SPOD modes at each frequency; $\{ \td{\Lambda}^{r_d}_i \}$, the $r_d$ SPOD energies at each frequency; $\{\td{G}_i \}$, the matrix with columns defined in \eqref{eq:gdef}; $\{ r_i \}$, the number of modes to be kept at each frequency; $\td{\Phi}$, the basis for reducing the initial condition and forcing terms.\\

{\textbf{Outputs:} $\{ \td{E}_i \}$, $\{ \td{H}_i \}$, $\{ \td{T}_i \}$, the operators in \eqref{eq:corr_pgp:final_eq} for each frequency; $\{ \td{C}^{\td{\Psi}}_i \}$, the operators that map the SPOD mode coefficients to $\ket{y}$ for each frequency.}

\end{algorithm}
\\
\begin{algorithm}
\caption{SSOP (online)}\label{alg:online}
\begin{algorithmic}[1]
\State \textbf{Input parameters:} $\ket{q}_0$,  ${\ket f}$, $\td{W}$, $\{ \td{\Psi}_i \}$, $\td{\Phi}$, $\{ \td{E}_i \}$, $\{ \td{H}_i \}$, $\{ \td{T}_i \}$, $\{ \td{C}^{\td{\Psi}}_i \}$
\State $\hat{\ket f} \gets \mathtt{FFT}({\ket f})$ \Comment{FFT of forcing} 
\State $\ket{a}^\td{\Phi}_0 = \td{\Phi}^* \td{W} \ket{q}_0$ \Comment{Reduced initial condition}
\State $\tilde{\ket{a}}^{\td{\Phi}}_0 \gets \td{0}_{p \times 1}$ \Comment{Initializing forcing sum term}
\For{$k \in \{1,2,\dots, N_{\omega} \}$}

\State $\ket{b}_k \gets \td{E}_k \hat{\ket f}_k$
\State $\tilde{\ket{a}}^{\td{\Phi}}_0  \gets \tilde{\ket{a}}^{\td{\Phi}}_0  + \frac{1}{N_\omega}\td{T}_k \ket{b}_k$ \Comment{Forcing sum}

\EndFor

\For{$k \in \{1,2,\dots, N_{\omega} \}$}

\State $\ket{a}_k \gets \ket{b}_k + \td{H}_k(\ket{a}^\td{\Phi}_0 - \tilde{\ket{a}}^{\td{\Phi}}_0 )$ \Comment{Assigning SPOD coefficients}
\State $\hat{\ket y}_k \gets \td{C}^{\td{\Psi}}_k \ket{a}_k$ \Comment{Constructing observable in frequency domain}
\EndFor

\State $\ket{y} \gets \mathtt{IFFT}(\hat{\ket y})$ \Comment{Observable in time domain}
\end{algorithmic} 
\vspace{1\baselineskip}
\textbf{Inputs:} $\ket{q}_0$, the initial condition; ${\ket f}$, the forcing as a function of time; $\td{W}$, the weight matrix; $\{ \td{\Psi}_i \}$, the retained SPOD modes for each frequency; $\td{\Phi}$, the basis for reducing the initial condition and forcing terms; $\{ \td{E}_i \}$, $\{ \td{H}_i \}$, $\{ \td{T}_i \}$, the operators in \eqref{eq:corr_pgp:final_eq} for each frequency; $\{ \td{C}^{\td{\Psi}}_i \}$, the operators that map the SPOD mode coefficients to $\ket{y}$ for each frequency. \\
\textbf{Output:} $\ket{y}$, the observable in the time domain.
\end{algorithm}
So long as the offline time is feasible, which we show next, the online scaling is the salient cost. In calculating the online cost, we count the operations necessary to go from the time domain forcing and initial condition to the SPOD coefficients. In practice, $\ket{y}$ is likely small, thus multiplying the SPOD coefficients by $\td{C} \td{\Psi}_k$ and taking the inverse DFT contributes insignificantly to the scaling. The complexity of the online algorithm is
\begin{equation} \label{eq:alg:scaling}
    \mathcal{O} ((r p + r N_f + N_f \log N_\omega ) N_\omega) \text{.}
\end{equation}
This time should be compared with the complexity of a POD-Galerkin model, which is $\mathcal{O}( (rN_f + r^2 )N_t)$. The two are similar, and the differences are due to how $p$ compares with $r$, how $N_t$ compares with $N_\omega$, and the constants involved. In our numerical experiments, we find that the SPOD method is slightly faster for equal rank (but much more accurate). In Appendix~\ref{App:DEIM} we detail a further approximation that removes the $N_f$ scaling from \eqref{eq:alg:scaling} by employing the discrete empirical interpolation method (DEIM) \citep{Chaturantabut10} to approximate the forcing, and other $N_x$-tall vectors, via sparse samplings. This can lead to significant speed-up in cases where $N_f$ is large.
\\

In Algorithm~\ref{alg:offline}, we have assumed that the SPOD modes for the $\hat{\ket{q}}$ and $\hat{\ket{g}}$ have already been obtained. These modes can be obtained using the techniques described in, e.g., Refs. \citep{Towne2018spectral,Schmidt19}, and the cost for this step is $\mathcal{O}(r_d^2N_xN_\omega)$, where again $r_d$ is the number of SPOD modes obtained from the data, which is the same as the number of temporal blocks formed in Welch's algorithm. If evaluation of $\td{A}\ket{q}$ scales linearly with $N_x$ (i.e., $\td{A}$ is sparse), then the offline scaling is $\mathcal{O}(r_d^2N_xN_\omega)$. If the evaluation of $\td{A}\ket{q}$ scales quadratically with $N_x$, then the offline scaling is $\mathcal{O}(r_d^2N_xN_\omega + r_dN_x^2N_\omega)$.
\\

\subsection{Data-free method} \label{sec:data_free}
Here, we briefly outline a data-free version of the method. Ref. \citep{Towne2018spectral} established a connection between SPOD modes, which come from data, and resolvent modes, which come directly from the system matrices. We first define the singular value decomposition
\begin{equation}
    \td{X}_q\td{R}_k \td{B} \td{X}_f^{-1} = \tilde{\td{U}}_k \td{\Sigma}_k \tilde{\td{V}}_k^* \text{.}
\end{equation}
Here, $\td{X}_q$ is the Cholesky factor of the weight matrix that defines the energy of the state, and $\td{X}_f$ is the Cholesky factor of the weight matrix that defines the norm of the forcing, i.e., $\td{W} = \td{X}_q^*\td{X}_q$, $\td{W}_f = \td{L}_f^*\td{L}_f$ and $\|\ket{f} \|^2 = \ket{f}^* \td{W}_f \ket{f}$. The resolvent response modes at frequency $\omega_k$ are then defined as $\td{U}_k = \td{X}_q^{-1} \tilde{\td{U}}_k$. The relation between these modes and the SPOD modes is the following: if the forcing at frequency $\omega_k$ is white in space with respect to the forcing norm, i.e., $\mathbb{E}[\hat{\vd{f}}_k \hat{\vd{f}}_k^*] = \alpha \td{W}_f^{-1}$ for some constant $\alpha$, then the SPOD modes are equivalent to the resolvent modes,
\begin{equation}
    \td{\Psi}_k = \td{U}_k  \text{,}
\end{equation}
and the SPOD energies are proportional to the square singular values,
\begin{equation}
    \td{\Lambda}_k = \alpha^2 \td{\Sigma}_k^2  \text{.}
\end{equation}
\\

Formally, this equivalence holds only if the forcing is white; however it has been demonstrated extensively in the fluid mechanics literature that the two sets of modes are often remarkably similar even when the forcing is far from white \citep{McKeon10,Beneddine16,Schmidt18,Nogueira19,Pickering20}. This can be leveraged to formulate a data-free ROM by substituting the resolvent modes for the SPOD modes at each frequency. We refer to this version of the method as resolvent SSOP.
\\

The scaling for a straightforward computation of the inverse that defines the resolvent, and the matrix exponential and inverse in \eqref{eq:corr_pgp_toapprox} is cubic in $N_x$. However, data-free methods \citep{Farghadan25} exist that can be used to drastically reduce this scaling, as is required for large systems. 
\\

The natural point of comparison for the resolvent SSOP method is balanced truncation, which is also data-free. We make this comparison in Section~\ref{sec:GL_results}, finding that at most parameter values, resolvent SSOP is more accurate. Balanced truncation shares the worst-case offline scaling, and approximations that reduce the scaling exist as well \citep{Willcox02, Rowley05}.

\section{Examples} 
\label{sec:results}
Here, we demonstrate the proposed method on two examples: a linearized Ginzburg-Landau problem and a scalar transport problem. The former is a dense system of dimension $N_x = 220$, and the latter is a sparse system of dimension $N_x = 9604$. For the Ginzburg-Landau system, we compare the accuracy and cost of the proposed method to those of POD-Galerkin projection and a statistics-enhanced version of balanced truncation. The error is roughly two orders of magnitude lower for the proposed method than the other two (depending on the case), and the CPU time is similar for all methods. For the scalar transport case, the offline time for balanced truncation makes the method (in its unapproximated form) infeasible, so we compare only to a POD-Galerkin model. The proposed method is again orders of magnitude more accurate at similar CPU cost.
\\

The SSOP method with an average of $r$ modes at each frequency is compared to the two time-domain methods with $r$ modes. In both time-domain methods, we use an integrator with an adaptive time step, but calculate the error at the temporal grid corresponding to the fastest frequency used in the SSOP method (which was longer than the time steps taken by the integrators in both examples). This means that the time-domain methods, as well as the proposed SSOP method, use $rN_\omega$ degrees of freedom to represent the solution over the interval. We will show that the SSOP and the two time-domain methods at the same $r$ are comparable in terms of CPU time, so the comparison that ultimately is the most meaningful is the accuracy of the methods at the same $r$.
\\

\subsection{Linearized Ginzburg-Landau problem} \label{sec:GL_results}
In continuous space, the complex linearized Ginzburg-Landau equation is
\begin{equation} \label{eq:GL_def}
    \dot{q}(x,t) =  \mathcal{A}q(x,t) + f(x,t) \text{,}
\end{equation}
where
\begin{equation}
    \mathcal{A} = -\nu \frac{\partial}{\partial x} + \gamma \frac{\partial^2}{\partial x^2} + \mu(x) \text{,}
\end{equation}
and $f(x,t)$ is a forcing. Following Ref. \citep{Bagheri09}, we set $\nu = 2 + 0.4i$ and $\gamma = 1 - i$. The parameter $\mu(x)$ takes the form
\begin{equation}
    \mu(x) = (\mu_0 - c_\mu^2) + \frac{\mu_2}{2}x^2 \text{,}
\end{equation}
with $c_\mu = 0.2$, $\mu_2 = -0.01$ \citep{Bagheri09}. The parameter $\mu_0$ is a bifurcation parameter; the linearized system transitions from global stability to global instability when $\mu_0$ exceeds $0.397$. We set $\mu_0 = 0.229$ \citep{Towne2018spectral} for the majority of our numerical experiments and later explore the effectiveness of the ROMs as $\mu_0$ varies from $0.079$ to $0.379$. The system can be interpreted as an advection-diffusion equation with a local exponential growth term. The equation supports traveling wave behavior in the positive $x$-direction and is stable in the sense that all the eigenvalues of the linear operator (discretized or continuous) are negative, so all solutions to the unforced equations decay asymptotically. Whether the exponential term promotes local growth or local decay depends on the sign of $\mu(x)$. With the parameters used, $\mu(x)$ is positive when $x \in [-6.15,6.15]$ and negative elsewhere, so as waves move through this region, they grow substantially before decaying once again after passing through it. 
\begin{figure}
    \centering
    \input{Figures/state_and_forcing_spacetime}
    \includegraphics{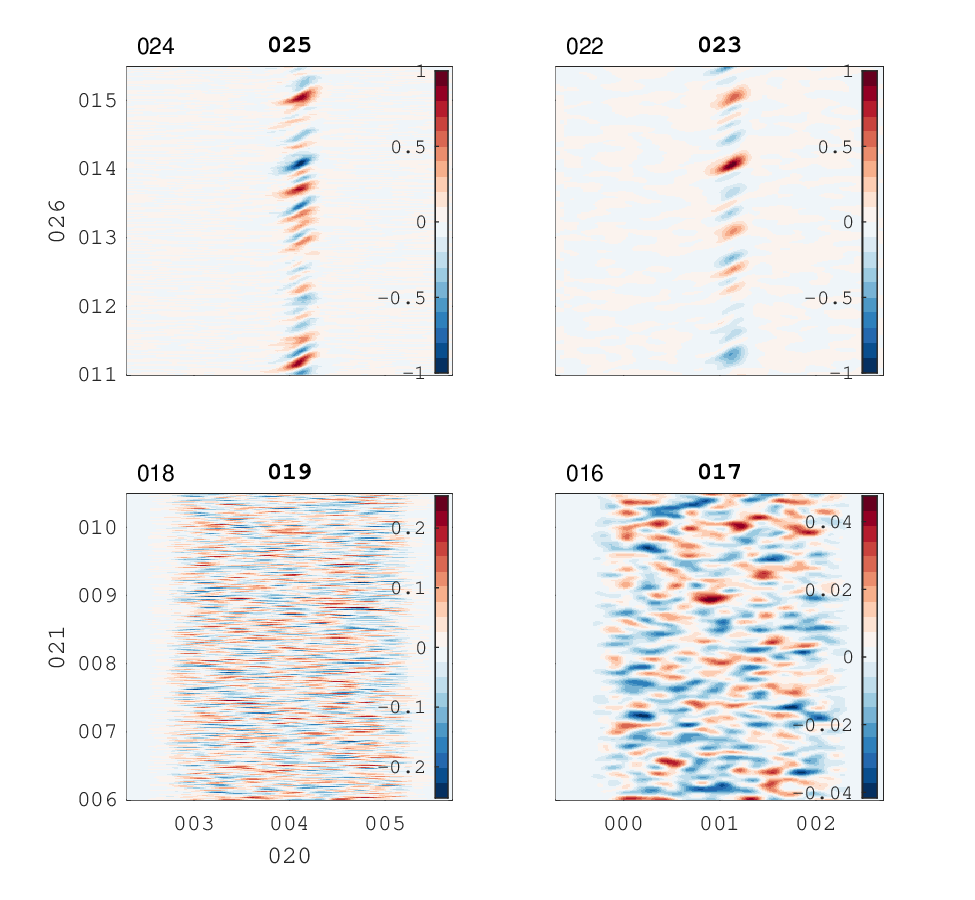}
    \caption{Ginzburg-Landau state and forcing trajectories: (a) the state resulting from the white forcing in (c); (b) the state resulting from the Gaussian forcing in (d). Both forcings have the same spatial correlation, but the short temporal correlation in the white forcing leads to more jagged structures in the corresponding state. Both trajectories consist of waves traveling in the positive $x$-direction that are amplified in a region near $x=0$.}
    \label{fig:state_and_forcing_spacetime}
\end{figure}
\\

When the equation is discretized in space, it takes the general form in (\ref{eq:ODE}). Following \citers{Bagheri09,Chen11}, we use a pseudo-spectral Hermite discretization with $N_x = 220$ collocation points \citep{Towne2018spectral}, and solve the discretized equations using MATLAB's \texttt{ode45}, an explicit Runge-Kutta (4,5) integrator. Following \citer{Towne2018spectral}, training data to compute the SPOD modes is generated from a single long run of the forced system comprising a short transient followed by $12000$ time steps of data with $\Delta t = 0.2$. The transient is discarded so that the modes are independent of the initial condition used in the run, and the $12000$ time steps are segmented into $r_d = 142$ overlapping blocks, each of length $N_\omega = 1024$ time steps. 
\\

Figure~\ref{fig:state_and_forcing_spacetime}(a,b) shows two space-time trajectories of the state $q$, each $1024$ time steps in length. The diagonally oriented structures demonstrate the traveling wave behavior of the system, and it is clear that the waves are amplified and then attenuated as they pass through $x=0$. These space-time trajectories are to space-time POD (and thus SPOD) as snapshots are to POD: the more structure there is in the trajectories, the fewer space-time modes are needed to accurately represent them. For example, in Figure~\ref{fig:state_and_forcing_spacetime}(a) the state is forced with band-limited temporally white noise and a Gaussian spatial correlation while the state in Figure~\ref{fig:state_and_forcing_spacetime}(b) is forced with a temporally, as well as spatially, Gaussian noise. The resulting state from the white forcing has more detailed structures -- a good proxy for higher rank behavior -- so trajectories with this forcing require more space-time modes to be accurately represented. The temporally white and temporally Gaussian forcings are shown in Figure~\ref{fig:state_and_forcing_spacetime}(c,d). Note that the forcing occupies the entire domain in this example, i.e., $N_f = N_x$. We thus choose $p = N_x$ vectors in the intermediary basis used to reduce the transient operator, because this will not have a large impact on the CPU time in this case. The SPOD method, therefore, scales like $N_x$, as do POD-Galerkin projection and balanced truncation, and there is no scaling benefit gained by using an intermediary basis, so we set it to the identity in this example. Despite this spatially extensive forcing, all three ROMs maintain a significant advantage over the FOM due to the latter being dense, resulting from the pseudo-spectral discretization.
\\
\begin{figure}[]
    \centering
    \input{Figures/romega3}
    \includegraphics{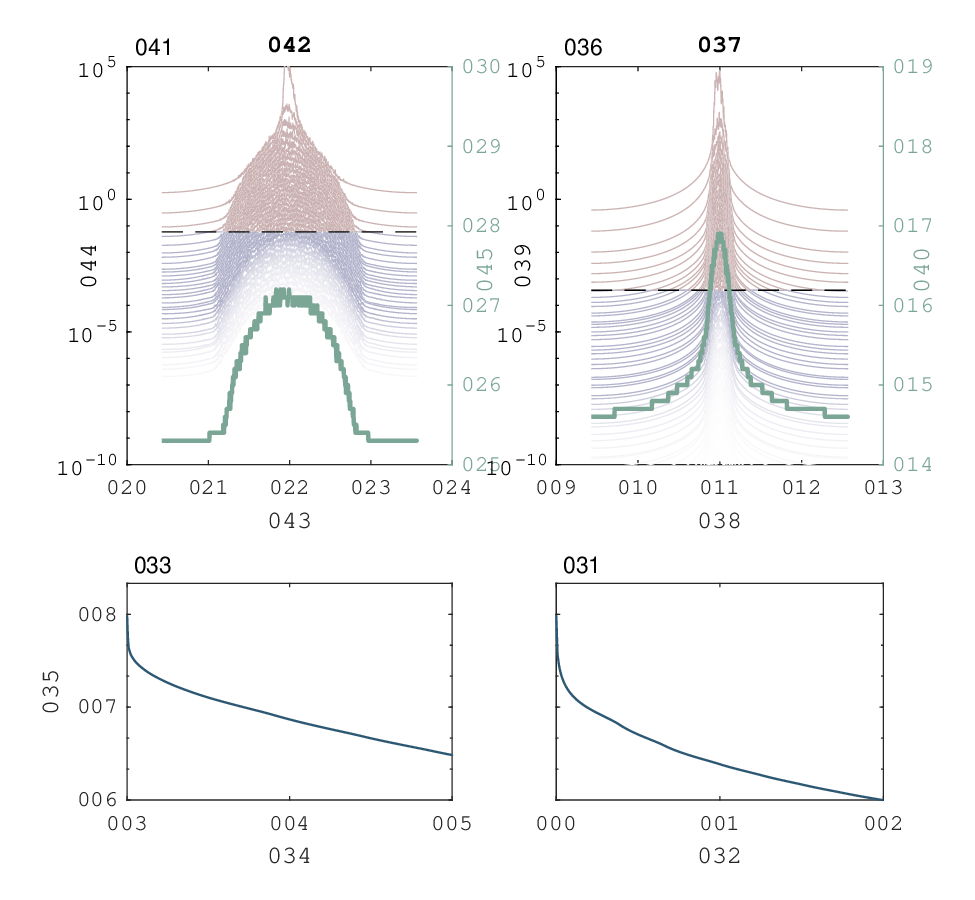}
    \caption{SPOD mode energies: (a-b, left axis) energy $\lambda$ of the retained and unretained modes. The top red curve is the energy of the first SPOD mode as a function of frequency $\omega$. The lower red and blue curves are the energies of the lower mode numbers, as functions of frequency. The retained modes (red) are the overall highest-energy modes, and the threshold (dashed) is determined as the energy of the $N_\omega r = 10240$-th most energetic mode; (a-b, right axis) number of modes that clear the threshold as a function of frequency. (c-d): the fraction of excluded energy (see \eqref{eq:excluded_energy}) as a function of $r$, the average number of modes per frequency.}
    \label{fig:romega}
\end{figure}

\begin{figure}[]
    \centering
    \input{Figures/ROM_state_error_snaps}
    \includegraphics{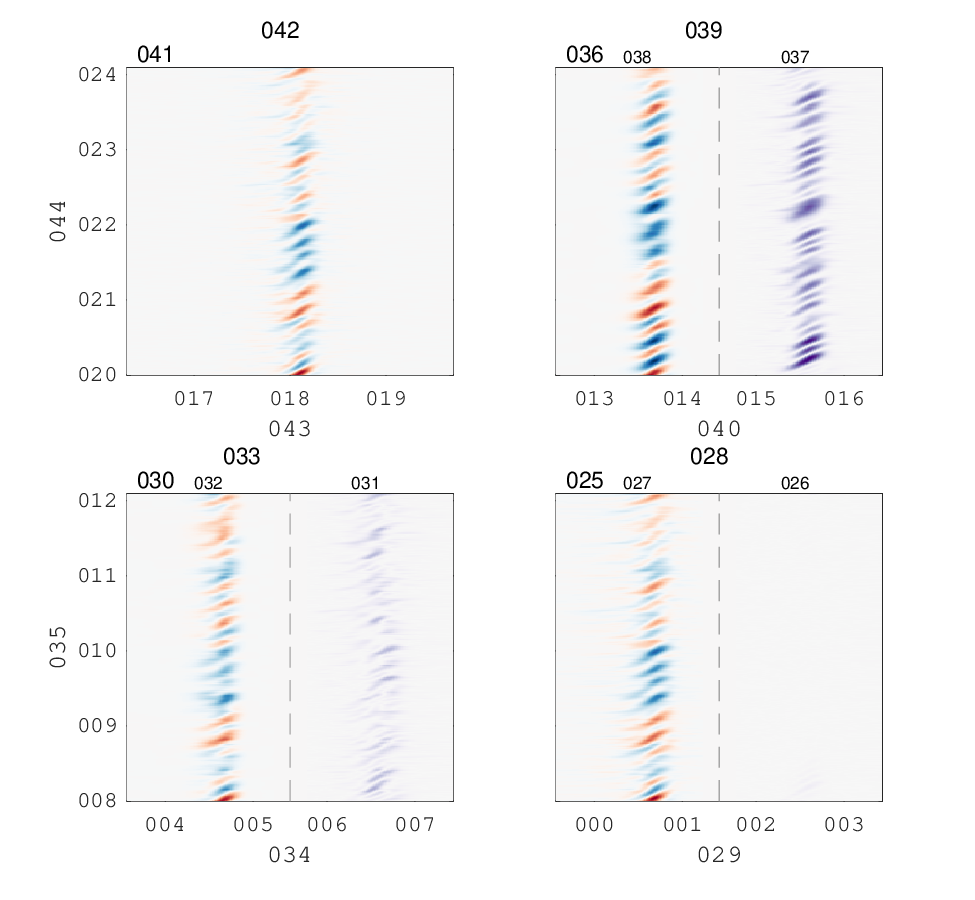}
    \caption{FOM trajectory (a) and ROM predictions thereof (b-d) along with the errors for the Ginzburg-Landau system. The error fields shown are the absolute value of the difference of the FOM and ROM trajectories. The peak error value (and the upper limit on the error color scale) is $87 \%$ of the peak absolute value of the state.}
    \label{fig:ROM_state_error_snaps}
\end{figure}

Figure~\ref{fig:romega} illustrates various features of the mode energies. The top panel shows the SPOD mode energies $\lambda$ as a function of $\omega$. Each curve is a particular mode number as a function of frequency. The decision to retain $N_\omega r$ total modes in the ROM selects an energy threshold below which modes are not retained (indicated by the dashed line in the figure). After ordering the energies of all mode numbers at all frequencies, the threshold is given by the $N_\omega r$-th energy $\tilde{\lambda}_{N_\omega r}$. At frequencies where a given mode number is above (red) the energy threshold (dashed), it is retained. Where it is below (blue) the energy threshold, it is truncated. The green curve shows the number of modes that meet or exceed this threshold as a function of $\omega$. For example, at the dominant frequency in the white noise case, $22$ modes are retained, whereas at the highest and lowest frequencies, only $3$ are retained. The bottom panel shows the fraction of energy that is excluded depending on the number of modes retained, as given by the formula
\begin{equation} \label{eq:excluded_energy}
    \frac{\sum_{k = 0}^{N_\omega - 1} \sum_{j = r_k + 1}^{N_x} \lambda_{kj}}{\sum_{k = 0}^{N_\omega - 1} \sum_{j = 1}^{N_x} \lambda_{kj}} \text{.}
\end{equation}
This quantity represents the fraction of the energy in the training data that is not representable by the retained modes. It depends on $r$, the average number of modes retained per frequency, and as $r$ is increased the excluded energy decreases as the SPOD basis at each frequency is enriched. For the Gaussian forcing case, the excluded energy fraction drops more quickly initially, indicating that the state is more accurately represented with a given number of SPOD modes in the Gaussian case relative to the white case.
\\

The test data comprise $173$ trajectories, again with $\Delta t = 0.2$ and with each trajectory of length $1024$ time steps. These test trajectories are generated as sub-trajectories from a long run, and the initial conditions and forcings used are not present in the training data. For each ROM, we calculate the error at the $1024$ points for each trajectory as the square norm of the difference with the FOM solution. Throughout this section, we compare the performance of the proposed SSOP method to that of POD-Galerkin projection and an enhanced version of balanced truncation. The reader is referred to \citer{Rowley17} for a description of POD-Galerkin projection. We use the MATLAB function \texttt{balreal} to generate the reduced system matrices and basis for balanced truncation, which is based on \citer{Moore81,Laub87}. In its usual form, balanced truncation does not make use of (i.e., does not require) data or knowledge of the statistics of the problem it is applied to. It may be improved with this information by `whitening' the forcing, i.e., transforming the system to one where the forcing is spatially white before performing the usual balanced truncation algorithm. In this application, where the forcing is far from spatially white, we observe that this variant of balanced truncation substantially outperforms the standard version, so we use it as a benchmark along with POD-Galerkin projection. We solve the reduced equations with the MATLAB function \texttt{ode45} for both methods, and run both the FOM and all ROMs using six cores of an Intel Xenon 6128 processor. The cost of building the ROM (with $r = 10$) was $0.7 \%$ of the cost generating the FOM data.
\\

Figure~\ref{fig:ROM_state_error_snaps} shows the three ROM approximations of a trajectory, along with the errors in these approximations. All ROMs here use $r = 2$ modes, and the error field shown here is the absolute value of the difference between the FOM and ROM trajectories, i.e., $| \tilde{q}(x,t) - q(x,t) |$ where $\tilde{q}(x,t)$ is the ROM result.. The enhanced balanced truncation produces a better result than POD-Galerkin projection, but the error of the SSOP method is much lower than both benchmarks. 
\\

Hereafter, we compute the error as a function of time as the square norm of the difference of the ROM and FOM solutions, averaged over the test trajectories and normalized by the mean square norm of the FOM solution, i.e., 
\begin{equation} \label{eq:error_t}
    e(t) = \frac{\sum_{i = 1}^{N_{\text{test}}} \| \tilde{\ket{q}}^i(t) - \ket{q}^i(t) \|_x^2 } { \frac{1}{T}  \sum_{i = 1}^{N_{\text{test}}} \int_0^T \| \ket{q}^i(t') \|_x^2 \ \text{d}t'} \text{.}
\end{equation}
We also report the mean of this quantity over time. 
\\

\begin{figure}[]
    \centering
    \input{Figures/Err_v_time_2ds}
    \includegraphics{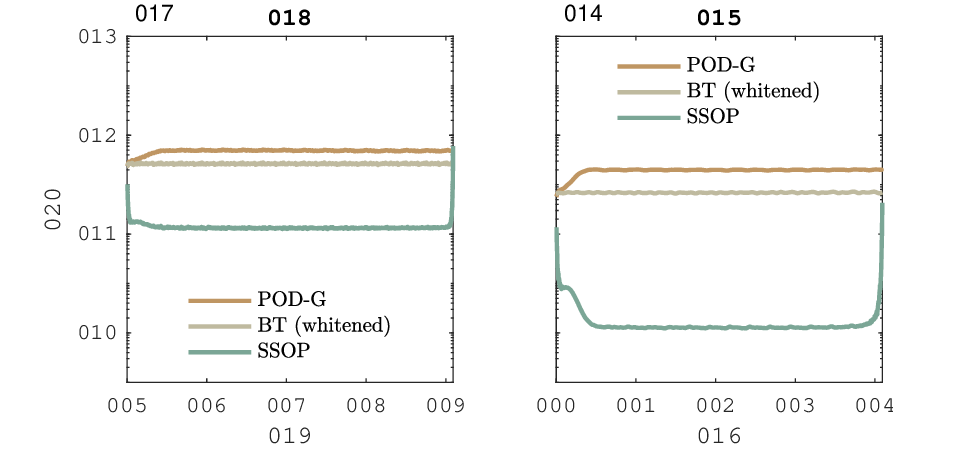}
    \caption{Error with $r = 10$ modes relative to the FOM for the Ginzburg-Landau system, averaged over $173$ trajectories. The large difference in accuracy is due, in large part, to the ability of the SPOD modes to represent trajectories more accurately than space-only modes. This difference is larger in the Gaussian forcing case.}
    \label{fig:err_v_time_2ds}
\end{figure}

Figure~\ref{fig:err_v_time_2ds} shows the error, defined in \eqref{eq:error_t}, for the various ROMs. For the POD-Galerkin and enhanced balanced truncation models, $10$ spatial modes are used for each of the $1024$ time steps. For the SSOP method, the $10 \times 1024 = 10240$ most energetic SPOD modes are used, and are distributed over the frequencies as shown in Figure~\ref{fig:romega}. Most notably, the error is nearly two orders of magnitude smaller for the SSOP method than it is for the other two methods. Given the analysis in the previous sections, this is not surprising: the SPOD modes are (nearly) optimal in that the representation error with some number of SPOD modes is smaller than (nearly) every other space-time basis. Again, this representation is recovered by the SSOP method up to the errors introduced from approximating the operators and the non-periodicity of the forcing on the temporal interval. The error from all methods is larger in the white-noise forcing case. This is to be expected because the resulting behavior of the state is higher rank in this case relative to the Gaussian forcing. The error of the SSOP method decreases by more in the Gaussian forcing case because it takes explicit advantage of the additional spatiotemporal coherence relative to the white forcing case. 
\\

\begin{figure}[]
    \centering
    \input{Figures/Err_v_modes_2ds}
    \includegraphics{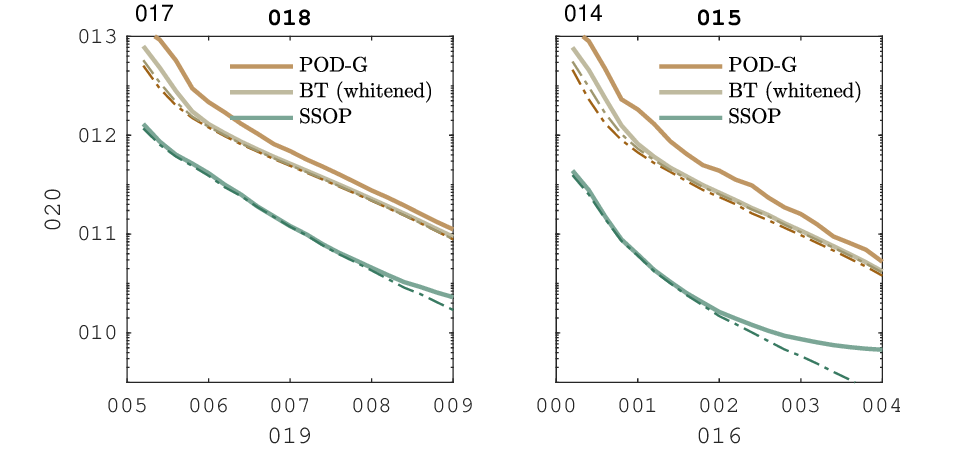}
    \caption{Error as a function of the number of modes for the Ginzburg-Landau system. The values reported here are the time averages of the error defined in \eqref{eq:error_t}.}
    \label{fig:err_v_modes}
\end{figure}

Next, we investigate the dependence of the error on the number of modes retained. Figure~\ref{fig:err_v_modes} shows, as one might expect, that the error in all methods decreases with the number of modes retained. The gulf between the SSOP method and the two space-only ROMs is roughly maintained over the range of modes shown for both the white and Gaussian forcing cases. The dashed lines are the projection of the full-order solution onto each respective set of modes, which we refer to as the representation error. For example, the dashed green line is the SPOD mode representation error, i.e., the error of the FOM solution projected into the span of the SPOD modes. We emphasize that the motivation for this work is the fact that the SPOD representation error is substantially below the POD representation error. The SSOP solution error and representation error are nearly identical before the accuracy of the former is limited by the full-order frequency domain error at around $16$ modes and $10^{-6}$ error. Until this point, the SSOP method indeed achieves the lowest error possible using SPOD modes, which is nearly the lowest error with any set of space-time modes. The SSOP solution error is not only lower than the POD-Galerkin and enhanced balanced truncation solution errors, but also the respective representation errors of these bases. The POD representation error is a lower bound for the error for any time-domain Petrov-Galerkin method, such as balanced truncation or least-squares Petrov-Galerkin \citep{Carlberg10}, because this is precisely the quantity that POD modes minimize. Indeed, the balanced truncation error is within this bound. We view the fact that the SSOP solution error is significantly below the POD representation error as one of the major achievements of this work.

\begin{figure}[]
    \centering
    \input{Figures/CPU_v_modes_2ds}
    \includegraphics{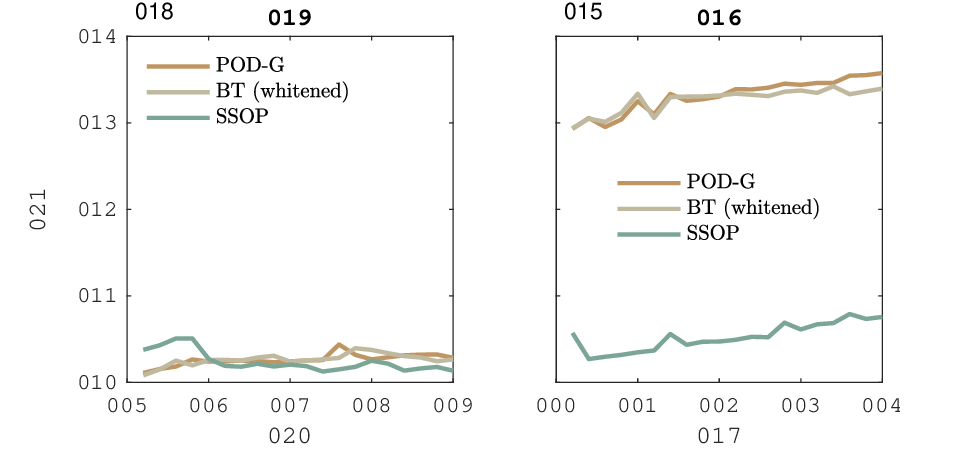}
    \caption{Average CPU time as a function of the number modes for the Ginzburg-Landau system. The values are normalized by the CPU time of the FOM.}
    \label{fig:cpu-v_modes_2ds}
\end{figure}

In Figure~\ref{fig:cpu-v_modes_2ds}, we show the CPU time as a function of the number of modes retained. All values are reported as a fraction of the FOM time of $3543$ seconds in total for the $173$ runs. The SSOP CPU time scales linearly with the number of modes retained, but here, there are too few modes to see this scaling. Nonetheless, the SSOP method is substantially faster than the two benchmarks and runs in roughly two thousandths of the FOM time. This time includes the time to take the Fourier transform of the forcing and the inverse Fourier transform of the response. The times for the white forcing cases and Gaussian forcings are comparable for the SSOP method, but the FOM takes substantially longer with the white forcing. 
\\

Thus far, we have tested the method on a system with the same statistics as the system on which the model was trained. In many applications, this assumption is too generous; a model trained on one system may be used to predict the behavior of a different system. We test the robustness of the method by training it on the Ginzburg-Landau system where the bifurcation parameter is $\mu_0 = 0.229$, then testing it on a range of Ginzburg-Landau systems with $\mu_0 \in [0.079,0.379]$ with increments of $0.03$. Physically, this range corresponds to shifting from modal behavior at the lower end to strongly non-modal behavior at the upper end \citep{Bagheri09}. The training system ($\mu_0 = 0.229$) is the center of this range, but its behavior is more similar to systems at the lower end than the higher end; one metric for this is the optimal transient growth \citep{Trefethen99,Schmid07}, which is just above $1$ for $\mu_0 = 0.079$, approximately $5$ for $\mu_0 = 0.229$, and nearly $200$ for $\mu_0 = 0.379$. We compare the performance of the SSOP model on the out-of-sample data to that of the POD-Galerkin model. For both, using training data from one system consists of simulating that system and obtaining the SPOD or POD modes from the data. To build the model on the test system, the ROM operators are computed with the modes from the training system and the $\td{A}$ and $\td{B}$ operators from the test system. We also compare results to the enhanced balanced truncation, which does not use data but does use the statistics of the forcing. 
\\

\begin{figure}[]
    \centering
    \input{Figures/vary_mu_dd}
    \includegraphics{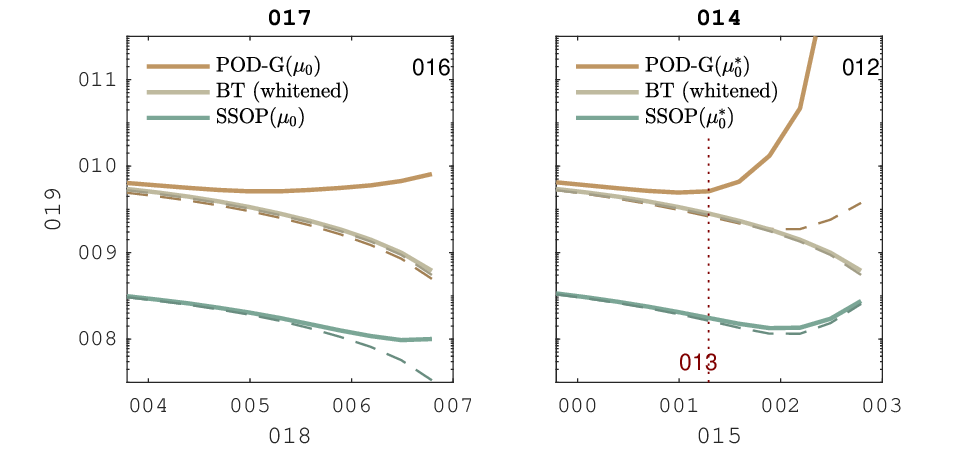}
    \caption{Data-driven ROM performance across a range of Ginzburg-Landau systems parameterized by $\mu_0$: (a) Models are trained and tested on the same system for a range of $\mu_0$; (b) Models are trained on a system with $\mu_0 = 0.229$ and tested on a range of $\mu_0$. All models use $r = 10$. The enhanced version of balanced truncation in this figure does not use data to obtain its test and trial bases, so the BT results are the same in (a) and (b). It is included among the data-driven results because it uses the forcing statistics to whiten the system. All values of $\mu_0$ shown are smaller than the critical value of $\mu_0 = 0.397$.}
    \label{fig:vary_mu_dd}
\end{figure}

Before considering out-of-sample conditions, in Figure~\ref{fig:vary_mu_dd}(a), we first show the performance where the ROMs are trained and tested on the same system, i.e., $\mu_0$ is the same. The error is computed with $r = 10$ for all models and is shown for the range of $\mu_0$. The dashed lines once again denote the projection of the FOM solution onto the respective modes; they are a lower bound for the respective methods. SSOP maintains its substantial accuracy superiority relative to the two baseline methods and to the POD projection over the range of $\mu_0$. Both the balanced truncation and SSOP models take advantage of the lower rank behavior at the high $\mu_0$ values, while the POD-Galerkin model is slightly less accurate due to the strong non-normality of $\td{A}$ at high $\mu_0$. 
\\

Figure~\ref{fig:vary_mu_dd}(b) shows the results for the same range of test systems, but where the POD-Galerkin and SSOP models use modes from $\mu_0 = 0.229$. The enhanced balanced truncation is the same between (a) and (b) because it only uses the statistics of the forcing, not the state, which are the same between the systems. Most notably, the SSOP model delivers low error over the range of test systems when it is trained on $\mu_0 = 0.229$. Even at high $\mu_0$, where the system is dominated by non-modal mechanisms that are far less prevalent in the training system \citep{Bagheri09}, the SSOP method is remarkably accurate. In fact, the SSOP method trained on the `wrong’ data ((b), solid green) is more accurate than the projection of the POD modes from the ‘right’ data ((a), dashed brown). That the system at, e.g., $\mu_0 = 0.379$ is substantially different from the one at $\mu_0 = 0.229$ may be seen by comparing the projection errors (dashed) for the POD and SPOD modes in (a) and (b) at $\mu_0 =  0.379$. For example, the error in the projection of the $\mu_0 = 0.379$ system onto the $\mu_0 = 0.379$ POD modes ((a), dashed brown) is an order of magnitude lower than the error in the projection of the $\mu_0 = 0.379$ system onto the $\mu_0 = 0.229$ POD modes ((b), dashed brown). The same can be said of the SPOD projection error between the two. The difference in the systems may also be seen in the fact that the POD-Galerkin model built from the $\mu_0 = 0.229$ POD modes is unstable at $\mu_0 = 0.379$ ((b), solid brown).
\\

Next, we investigate the robustness of the model to out-of-sample forcings using the following two tests. First, we train the model with $r = 5$ on the data generated by the white forcing and test it on the Gaussian forcing. We compare the statistics of the model output over $142$ blocks, namely the variance as a function of $x$ and the power spectral density as a function of $\omega$. Second, we use the same white-forcing-trained model to predict trajectories for four different non-stationary forcings. Unlike the white and Gaussian forcings, which are active in the entire domain, these non-stationary forcings are only nonzero in the interval $x\in [-12,12]$. Within this region, the forcings are given by $f(x,t) = T(t)X(x)$ where $X(x)$ is a Gaussian bump centered at $x = -10$. We test four choices of $T(t)$. 
\\

\begin{figure}[]
    \centering
    \input{Figures/psd_tke}
    \includegraphics{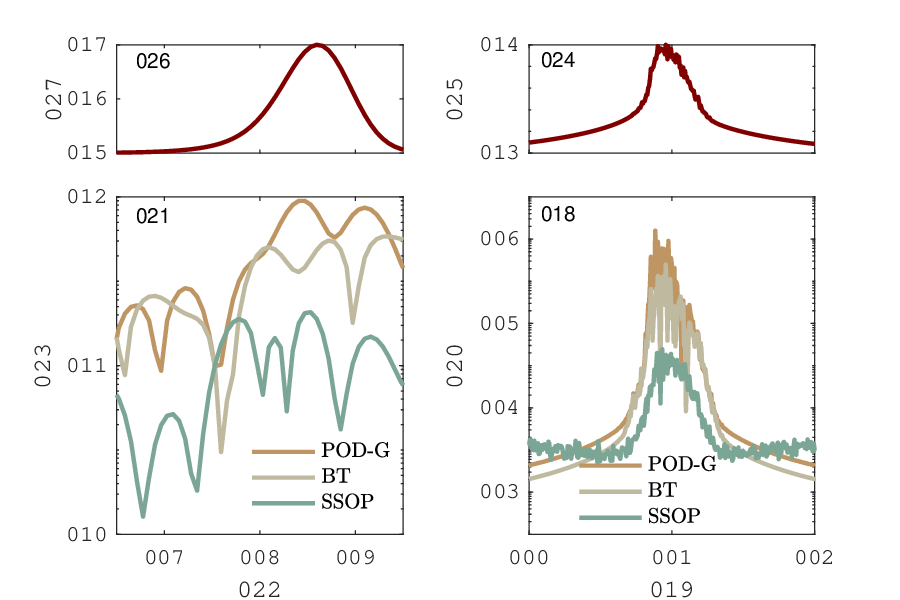}
    \caption{Error in the prediction of statistics for the Ginzburg-Landau system: (a) the FOM variance (TKE) as a function of $x$; (b) the power spectral density (PSD); (c) the ROM errors in predictions of the TKE; (d) the ROM errors in predictions of the PSD. The ROMs are built using data from the white-noise forcing.}
    \label{fig:psd_tke}
\end{figure}

Figure~\ref{fig:psd_tke} shows the results of the first test. The FOM variance, i.e., the turbulent kinetic energy, is shown in panel (a) and the errors in the ROM predictions thereof are shown in panel (c). Likewise, the FOM power spectral density is shown in panel (b) and the errors in the ROM PSDs are shown in panel (d). Both the mean energy and the PSD are normalized such that their peak is $1$, and both errors are computed as the absolute value of the difference between the FOM and ROM statistics. For both quantities, the proposed method produces more accurate statistics for the out-of sample forcings than does either POD-Galerkin projection or balanced truncation. We also note that the error in the SSOP prediction of the state for this test is $84$ and $20$ times lower than the POD-Galerkin and balanced truncation predictions, respectively. 
\\

\begin{figure}[]
    \centering
    \input{Figures/diff_forcings}
    \includegraphics{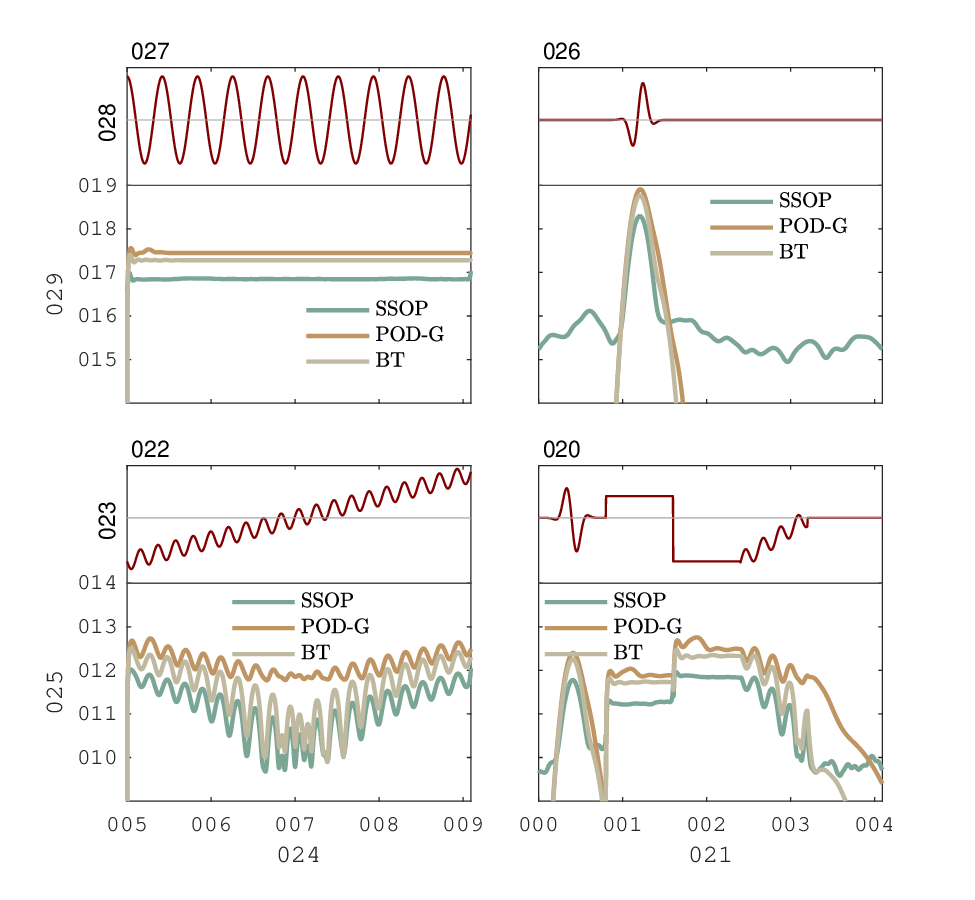}
    \caption{Results of the models trained with the white forcing on four out-of-sample non-stationary forcings to the Ginzburg-Landau system. The time-dependence of these forcings is shown in the top panel, and the resulting error for the three ROMs in the bottom panel. Respectively, the mean errors for the SSOP, POD-Galerkin, and balanced truncation models are (a) $.007$, $.03$, $.02$; (b) $.008$, $.03$, $.02$; (c) $.002$, $.02$, $.007$; (d) $.002$, $.01$, $.007$.}
    \label{fig:diff_forcings}
\end{figure}

Figure~\ref{fig:diff_forcings} shows the four non-stationary forcings, along with the error in the ROM predictions of the state for each. In each case, the errors are normalized by the mean square energy of the solution over the interval. The SSOP error is lower than that of the other two methods for all tests, but the difference is smaller than in previous tests. This likely may be attributed to the fact that the support of the forcing is quite different than in the training data for the model. Nevertheless, the method gives relatively accurate predictions for this out-of-sample forcing.

\begin{figure}[]
    \centering
    \input{Figures/vary_mu_df}
    \includegraphics{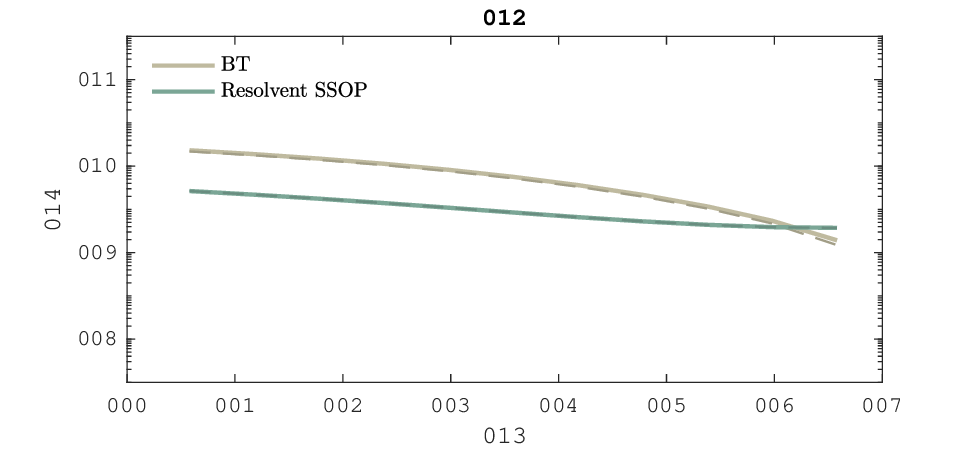}
    \caption{The accuracy of the resolvent SSOP model and unwhitened balanced truncation for Ginzburg-Landau systems with a range of $\mu_0$. Both methods are data-free; the resolvent modes come directly from the system matrices in the governing equation, and ‘unwhitened’ indicates that we did not transform the system to one where the forcing is white using forcing data before performing balanced truncation. The dashed lines are the projections of the FOM solution onto the respective modes.}
    \label{fig:vary_mu_df}
\end{figure}

Figure~\ref{fig:vary_mu_df} compares the accuracy of the resolvent SSOP model, described in Section~\ref{sec:data_free}, to that of non-whitened balanced truncation, both with $r=10$. The test is conducted for the same range of $\mu_0$ values as used in \ref{fig:vary_mu_dd}. The non-whitened version of balanced truncation does not transform the system to one where the forcing is spatially white, and, as such, is data-free. Notably, the resolvent SSOP method, which is also data-free, outperforms balanced truncation, the current gold-standard for linear model reduction, over most of the range. The online algorithm for resolvent SSOP is identical to SSOP, so its cost is also comparable to balanced truncation. We also note resolvent SSOP nearly recovers the projection of the FOM solution onto the resolvent modes, as can be seen by comparing the solid and dashed green lines in Figure \ref{fig:vary_mu_df}.

\begin{figure}[]
    \centering
    \input{Figures/ST_snaps}
    \includegraphics{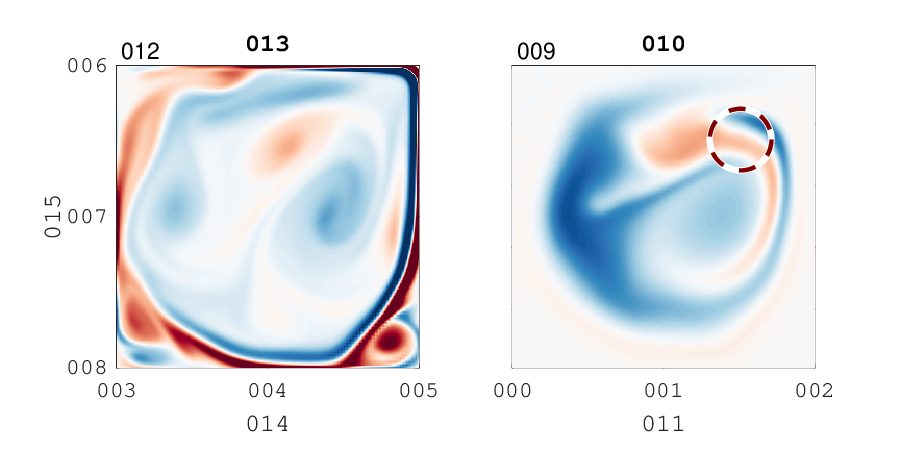}
    \caption{Visualization of the scalar transport problem: (a) the vorticity of the steady velocity field that transports the scalar; (b) a snapshot of the scalar field with the forcing region highlighted.}
    \label{fig:ST_snaps}
\end{figure}

\begin{figure}[!p]
    \centering
    \input{Figures/ST_error_fields}
    \includegraphics{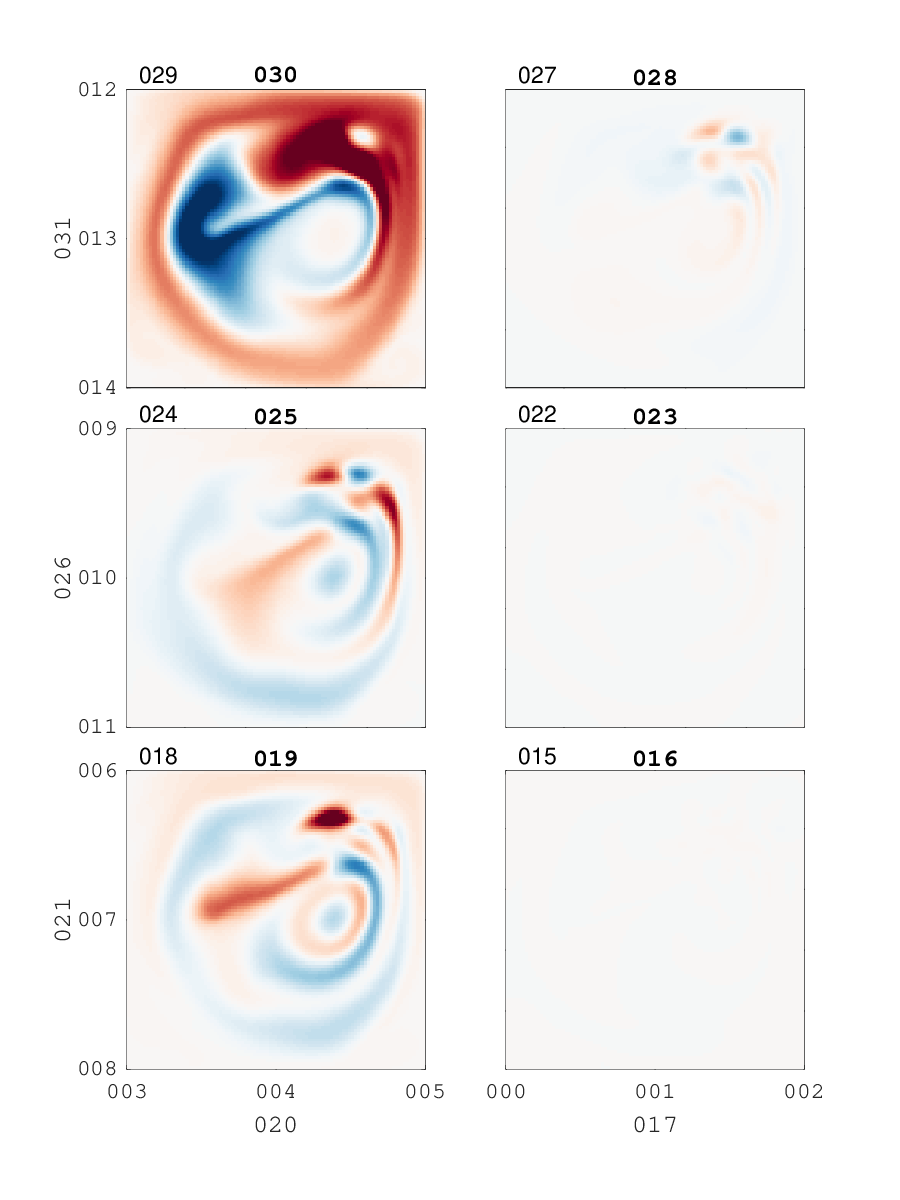}
    \caption{Error snapshots for the scalar transport problem: (a, c, e) POD-Galerkin ROM error for $2$, $10$, and $20$ modes; (b, d, f) SSOP ROM error for the same numbers of modes. The extremes of the color scale are $31 \%$ of the maximum absolute value of the FOM solution, and the maximum in the $2$-mode POD-Galerkin error field is $62\%$ of the maximum absolute value of the FOM solution.}
    \label{fig:ST_error_fields}
\end{figure}

\subsection{Scalar transport problem}
Next, we demonstrate the proposed algorithm on an advection-diffusion system modeling the transport of a scalar quantity in a steady fluid flow. The flow profile is the mean of a lid-driven cavity flow simulation at $Re = 30,000$. This problem differs from the Ginzburg-Landau example in three important ways: it is substantially larger ($N_x = 9604$ as opposed to $220$ in the Ginzburg-Landau case), the matrix $\td{A}$ is sparse, and the forcing occupies only a subset of the domain. The former means that the model is too large to compute the matrix operations without the approximations we described earlier. With this large $N_x$, computing the Gramians in balanced truncation is too costly as well, so we do not compare to it here, though we do note that there are effective data-driven approximations of it as well \citep{Willcox02,Rowley05}. A balanced truncation model must have greater error than the POD representation error, which we do report, and can be expected to share the CPU time of a POD-Galerkin model. That the forcing does not occupy the entire domain means that there is a CPU time savings in using the intermediary basis described in Section \ref{sec:transient_op}. 
\\

The continuous governing equations for the scalar transport case may be written in the same form as (\ref{eq:GL_def}), where $\mathcal{A}$ is now defined as
\begin{equation}
    \mathcal{A} = -\ket{u}(\td{x}) \cdot \nabla  + \eta \nabla^2  \text{,}
\end{equation}
and where $\ket{u}(\td{x})$ is the mean flow in the lid-driven cavity. We take $\eta = 0.001$ and the velocity of the lid to be Gaussian in $x$ (as opposed to a constant) to avoid discontinuities at the upper corners. The problem is nondimensionalized such that the maximum speed, occurring in the center of the lid, is $1$. We prescribe a forcing that is stochastic with Gaussian spatial and temporal autocorrelation in a region of Gaussian support centered at $\overline{\ket{x}} = [0.75,0.25]^T$; its statistics are given by 
\begin{equation} \label{eq:ST forcing stats}
    C_{ff}(\ket{x}_1 \ket{x}_2,t_1,t_2) = \exp \Bigg[- \left( \frac{|\ket{x}_1 - \overline{\ket{x}}|^2 + |\ket{x}_2 - \overline{\ket{x}}|^2}{l^2} + \frac{|\ket{x}_2 - \ket{x}_1|^2}{\xi^2} + \frac{(t_2 - t_1)^2}{\tau^2} \right) \Bigg] \text{,}
\end{equation}
where $| \ \cdot \ |$ is Euclidean distance. Here, $l = 0.1$ is the spatial width of the support of the forcing, $\xi = 0.07$ is the spatial correlation length, and $\tau = 1$ is the temporal correlation length. We use a second-order finite difference discretization with $98$ points in both directions and Dirichlet boundary conditions so as to mimic a heat bath. The FOM is solved using MATLAB's \texttt{ode45}. The vorticity of the underlying velocity field and a snapshot of the transported scalar are shown in Figure~\ref{fig:ST_snaps}. The red dashed circle in the latter indicates the region in which the equations are substantially forced -- it is the radius at which the autocorrelation of the forcing defined by \eqref{eq:ST forcing stats} drops by a factor of $e$ from its peak value. 
\\

For this example, we gather $50,000$ time steps of FOM data spaced $\Delta t = 0.5$ apart from which to calculate SPOD modes and approximate operators. Time is nondimensionalized such that it takes one time unit for the lid to cross the cavity. In this example, the time step used in the FOM integration was much smaller -- more than an order of magnitude for most times -- due to the stiffness of the system. The data was then segmented into $648$ overlapping blocks, each $256$ time steps in length, so $N_\omega = 256$, and $T = 128$. The test data is $128$ trajectories of the system, and error is defined in the same way as before -- as the square norm of the difference of the ROM and FOM solutions averaged over all test trajectories and normalized by the mean square norm of the solution itself. Both the ROMs and the FOM use 56 cores on a pair of Intel Xeon 6242R processors. The cost of building the $10$-mode ROM, including $168$ seconds to compute the SPOD modes, was $30 \%$ of the cost of generating the FOM data. The additional offline time for the POD-Galerkin ROM was trivial in comparison to the cost of running the FOM to generate data from which to obtain the POD modes. Finally, we use $p = 50$ for the intermediary basis.
\\

Figure~\ref{fig:ST_error_fields} shows error field snapshots for the two ROMs for $2$, $10$, and $20$ modes. The error field for the POD-Galerkin model with $20$ modes is larger than that for the proposed model with $2$ modes.
\\

\begin{figure}[!h]
    \centering
    \input{Figures/st_nodeim_error}
    \includegraphics{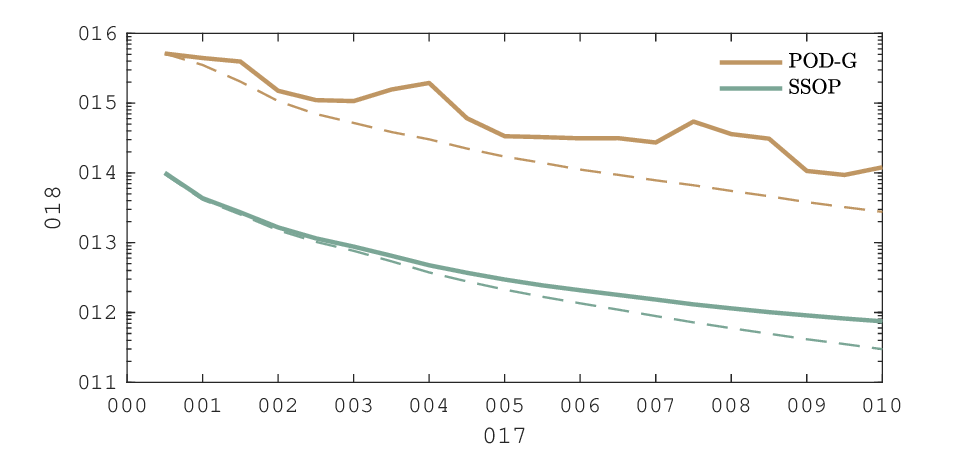}
    \caption{Accuracy of the proposed method compared to that of a POD-Galerkin model applied to the scalar transport problem. The dashed lines are the error of the full-order solution projected on the respective bases and are lower bounds for the error of the methods.}
    \label{fig:st_err_vs_modes}
\end{figure}

Figure~\ref{fig:st_err_vs_modes} shows the accuracy over a range of mode numbers. Once again, the proposed method produces error two orders of magnitude lower than that of the POD-Galerkin model. The dashed lines are again the representation error. The error of the POD projection is a lower bound for the POD-Galerkin error, as well as any spatial Petrov-Galerkin method, such as balanced truncation. We see that in this problem, like in the Ginzburg-Landau problem, the SSOP model far outperforms even this bound and is indeed fairly close to its own error bound.

\begin{figure}[!h]
    \centering
    \input{Figures/st_deim_timing}
    \includegraphics{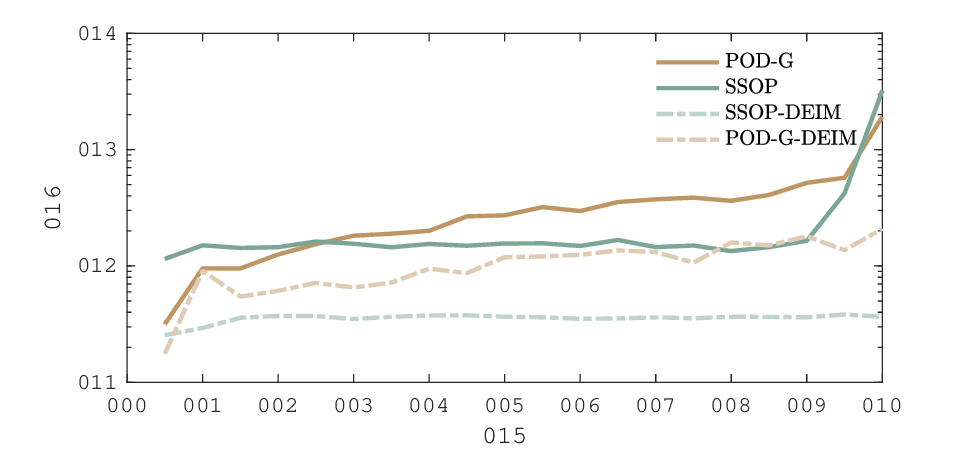}
    \caption{Average CPU time to solve the scalar transport problem as a function of the number of space-time modes used as a fraction of the FOM CPU time. The timing results from using DEIM to remove the $N_f$ scaling from the complexity of the model (and from the POD-Galerkin model) are also shown, and the details for this method are given in Appendix~\ref{App:DEIM}.}
    \label{fig:st_deim_timing}
\end{figure}

Figure~\ref{fig:st_deim_timing} shows the CPU time required to solve the scalar transport problem. The SSOP model is slightly faster in this case than the POD-Galerkin one. Too few modes are used to see the asymptotic linear scaling with the number of modes. Timing results from using DEIM to remove the $N_f$ scaling are also shown for both methods, and the DEIM-augmented version of the method is described in Appendix~\ref{App:DEIM}. There, we show that the error is not meaningfully affected by the DEIM approximation, but we also discuss some drawbacks of the augmented version of the method.

\section{Conclusions} \label{sec:Conclusion}
Space-time bases allow for a more accurate representation of a trajectory than do space-only bases with the same number of coefficients. In particular, the SPOD encoding of a trajectory with some number of coefficients may be orders of magnitude more accurate than the POD encoding of the same trajectory with the same number of coefficients. The obvious objectives are, therefore, to solve for these coefficients quickly and accurately, and we have pursued these for linear time-invariant systems.
\\

The method works as follows. The SPOD coefficients $\ket{a}_k$ at frequency $\omega_k$ are given by the Fourier transform of the trajectory $\hat{\ket{q}}_k$ at that frequency left-multiplied by the (weighted) transpose of the SPOD modes. In a linear system, $\hat{\ket{q}}_k$ obeys a linear relation involving the forcing (at all frequencies) and the initial condition. We derive this relation and left-multiply it by the transpose of the SPOD modes, which amounts to a projection of a solution operator onto SPOD modes. The resulting matrices are small and may be precomputed, and the online phase of the method involves small matrix-vector products to obtain the SPOD coefficients given the forcing and initial condition.
\\

We show, via two examples, that the SSOP method can indeed accomplish both objectives outlined above: it takes comparable CPU time to benchmark time-domain methods like POD-Galerkin projection and balanced truncation, and the solution from SSOP is roughly two orders of magnitude more accurate than both benchmarks. In fact, the SSOP solution is nearly two orders of magnitude more accurate than the projection of the FOM solution onto the POD modes, which is the lower bound on error for any time-domain Petrov-Galerkin method. We also demonstrated the robustness of the method by training and testing it on different systems. 
\\

A few negative aspects of the method are worth mentioning. The most limiting is that the method requires the entire forcing over the time interval of interest to be known before beginning the computation; the first step of the method is to take a FFT of the forcing, which cannot be done without the entire forcing in time. For some applications, this prevents the method from being applicable, while for others, it is not a problem. Second, the method works on a preprescribed interval $[0,T]$. If one wishes to obtain the solution longer than this interval, one can repeat the method with the value at the end of the interval as the initial condition. This is more cumbersome than extending the solution in a time-domain method. Finally, SPOD modes require more training data than POD modes, which limits the applicability of the proposed method in cases where training data is scarce. Where these disadvantages are not obstacles, however, we have shown the proposed method to be substantially more accurate at the same CPU cost compared to standard methods for linear model reduction. We hope that this will aid in applications of linear model reduction and increase interest in space-time methods. 
\\

Finally, we also described a resolvent-based version of the SSOP method that is data-free. Leveraging the ability of resolvent modes obtained directly from the linear system to approximate SPOD modes \citep{Towne2018spectral}, the resolvent SSOP method uses resolvent modes in place of SPOD modes within the SSOP method, eliminating the need for training data. We showed this method to be more accurate than balanced truncation, the state of the art for data-free model reduction. We see this result as quite promising, calling for further work comparing the method to balanced truncation. We also believe that this result underscores the promise of space-time techniques for model reduction.
\\

One clear area for further inquiry is whether the method can be generalized to the case of nonlinear governing equations. The method uses the fundamental solution to the LTI system, and, of course, no analogous solution exists for general nonlinear systems. One possible remedy is to treat the nonlinearity that would result from the present SPOD coefficients as an additional forcing on the system, and then again use the fundamental LTI solution. Nonlinearity as a forcing to a linear system is the perspective taken by resolvent analysis \citep{McKeon10}, and the successes of that field in predicting turbulent structures are cause for optimism that this perspective would be useful for model reduction. With such an approach, \eqref{eq:corr_pgp:final_eq} would become a coupled system of nonlinear algebraic equations, which would be solved online for the SPOD coefficients given the initial condition and forcing.
\section{Acknowledgments}
P.F. and A.T. gratefully acknowledge funding from the National Science Foundation grant no. 2237537. C.L. and O.S. gratefully acknowledge funding from the National Science Foundation grant no. 2046311.

\begin{appendices}


\section{Fourier representation of the time derivative} \label{App:Derivative_coefficient}
Here, we illustrate why substituting $i\omega_k \hat{\ket{q}}_k$ for $\dot{\ket{q}}$ when going to the frequency domain is only correct if $\ket{q}$ starts and ends at the same value on the interval. It is easier to do this using the continuous definition of the Fourier coefficients, rather than the DFT, so we compute the integral
\begin{equation}
    \hat{\dot{\ket{q}}}_k = \fint \dot{\ket q} e^{-i \omega_k t} \ dt\text{.}
\end{equation}
This integral may be solved by parts,
\begin{equation}
    \hat{\dot{\ket{q}}}_k = \frac{1}{T} \ket{q} \big|_{0}^{T} - \fint -i \omega_k \ket{q}(t) e^{-i \omega_k t} dt  \text{.}
\end{equation}
The boundary term gives $\frac{\Delta \ket{q}}{T}$ \cite{Martini21}, where $\Delta \ket{q} = \ket{q}(T) - \ket{q}(0)$, and the integral is the naive term $i\omega_k \hat{\ket q}_k$. The Fourier representation of the derivative is thus
\begin{equation}
    \hat{\dot{\ket{q}}}_k = i\omega_k \hat{\ket q}_k + \frac{\Delta \ket{q}}{T} \text{.}
\end{equation}

\section{DEIM-augmented algorithm} \label{App:DEIM}
Here, we present a means of improving the cost scaling in cases where the dimension of the forcing $N_f$ is large but the forcing is spatially structured. The idea is to use a sparse sampling of the forcing vectors, which are size $N_f$, using the discrete empirical interpolation method (DEIM) \citep{Chaturantabut10}. We also sparse sample the two vectors that are size $N_x$, the initial condition and the forcing sum terms in \eqref{eq:simpler_pgp_noib}, though these terms can also be handled using an intermediary basis, as before. The rank-$p$ DEIM approximation of a vector $\ket{v} \in \mathbb{C}^{N_x}$ is $\ket{v} \approx \td{U}_{\ket{v}} (\td{P}^T_{\ket{v}}\td{U}_{\ket{v}})^{-1}\td{P}^T_{\ket{v}} \ket{v}$, where the columns of $\td{U}_{\ket{v}} \in \mathbb{C}^{N_x \times p}$ are the POD modes for the ensemble from which  $\ket{v}$ is a sample, and $\td{P}^T_{\ket{v}} \in \{0,1\}^{p \times N_x}$ samples $p$ elements from $\ket{v}$ and is formed via the DEIM algorithm. 
\\

The DEIM algorithm is run for the forcing in the time domain, giving a set of sample points $\td{P}^T_{\ket{f}} \in \mathbb{C}^{p \times N_f}$ from which the forcing can be reconstructed accurately. The structures in the forcing at different frequencies will, in general, be different, so it is best to use a different spatial basis at each frequency to complete the DEIM approximation. We label the spatial basis for the forcing at the $k$-th frequency $\td{U}_{\hat{\ket{f}}_k}$. The approximation for the $k$-th forcing is $\hat{\ket{f}}_k \approx \td{U}_{\hat{\ket{f}}_k}(\td{P}^T_{\ket{f}}\td{U}_{\hat{\ket{f}}_k})^{-1} \td{P}^T_{\ket{f}} \hat{\ket{f}}_k$, so the first term in \eqref{eq:simpler_pgp_noib} is now 
\begin{equation}
    \td{E}_k \hat{\ket{f}}_k \approx \td{K}_{\hat{\ket{f}}_k}  \td{P}^T_{\ket{f}} \hat{\ket{f}}_k \text{,}
\end{equation}
where the matrix $\td{K}_{\hat{\ket{f}}_k} = \td{E}_k \td{U}_{\hat{\ket{f}}_k}(\td{P}^T_{\ket{f}}\td{U}_{\hat{\ket{f}}_k})^{-1} \in \mathbb{C}^{r_k \times p}$ is precomputed for each frequency. Note that the same sampling is used for each frequency; if the sampling were different for each frequency, one would need to take the DFT of the entire forcing (or the union of all the samplings), which would negate the scaling benefit of sparse sampling. This approach for the forcing term can be used in conjunction with the intermediary basis approach for the initial condition and forcing sum terms in \eqref{eq:simpler_pgp_noib}. These latter terms can also be handled with DEIM, which we show now.  
\\

The initial condition and forcing sum terms can be approximated with DEIM matrices $\td{U}_{\ket{q}_0}$ and $\td{P}^T_{\ket{q}_0}$ for the initial condition, and  $\td{U}_{fs}$ and $\td{P}^T_{fs}$ for the forcing sum. For the former, these matrices are formed by gathering all initial conditions within the training data and running the DEIM algorithm to obtain $\td{U}_{\ket{q}_0}$ and $\td{P}^T_{\ket{q}_0}$. From these matrices, the initial condition multiplied by $\td{F}_k$ is approximated as 
\begin{equation}
    \td{F}_k \ket{q}_0 \approx \td{K}_{\ket{q}_0} \td{P}^T_{\ket{q}_0} \ket{q}_0 \text{,}
\end{equation}
where  $\td{K}_{\ket{q}_0,k} = \td{F}_k \td{U}_{\ket{q}_0} (\td{P}^T_{\ket{q}_0} \td{U}_{\ket{q}_0})^{-1} \in \mathbb{C}^{r_k \times p}$ is precomputed. Similarly, $\td{U}_{fs}$ and $\td{P}^T_{fs}$ are obtained by running DEIM on the set of forcing sums, which must be calculated from each trajectory in the training data. With these matrices, the forcing sum multiplied by $\td{F}_k$ is approximated as
\begin{equation}
    \td{F}_k \frac{1}{N_\omega} \sum_{l} \td{\Psi}_l \td{E}_l \hat{\ket{f}}_l \approx \frac{1}{N_t} \td{K}_{fs,l} \sum_{l}  \td{T}_{l,fs} \td{E}_l \hat{\ket{f}}_l \text{,}
\end{equation}
where $\td{K}_{fs,l} = \td{F}_k \td{U}_{fs} (\td{P}^T_{fs} \td{U}_{fs})^{-1} \in \mathbb{C}^{r_k \times p}$ and $\td{T}_{l,fs} = \td{P}^T_{fs}\td{\Psi}_l \in \mathbb{C}^{p \times r_l}$ are both precomputed. These operators are kept separate (as opposed to multiplied as a precomputation step) to avoid $N_\omega^2$ scaling. With these approximations, the DEIM-augmented equation is
\begin{equation}
    \ket{a}_k =  \td{K}_{\hat{\ket{f}}_k}  \td{P}^T_{\ket{f}} \hat{\ket{f}}_k + \td{K}_{\ket{q}_0} \td{P}^T_{\ket{q}_0} \ket{q}_0 - \frac{1}{N_t} \td{K}_{fs,l} \sum_{l}  \td{T}_{l,fs} \td{E}_l \hat{\ket{f}}_l \text{.}
\end{equation}
This method removes the $N_f$ scaling and replaces it with $p$, thus, the DEIM-augmented version of the algorithm scales like $\mathcal{O}( (rp + p\log N_\omega)N_\omega )$. 
\\
\begin{figure}[!hbt]
    \centering
    \input{Figures/st_deim_error}
    \includegraphics{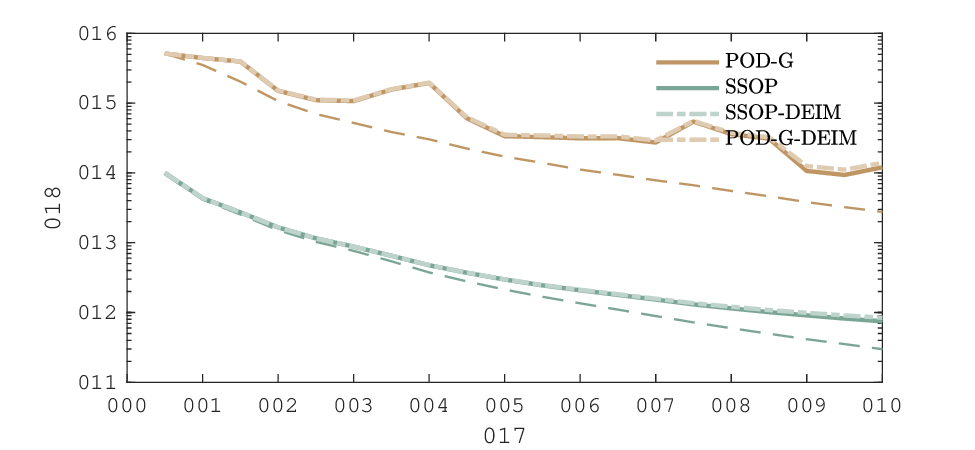}
    \caption{Accuracy of the DEIM-augmented version of the method compared to POD-Galerkin projection (and the DEIM version thereof) applied to the scalar transport problem. As long as enough sample points are used, DEIM does not introduce additional error. Again, the dashed lines are the error of the full-order solution projected on the respective bases.}
    \label{fig:st_deim_error}
\end{figure}
\\

To demonstrate the augmented version of the method, we apply it to the scalar transport problem described in the main text. The forcing in this problem meets the conditions for a large improvement: $N_f = 2050$ is large, but the forcing is spatially structured, coming from \eqref{eq:ST forcing stats}, so it can be approximated via sparse sampling effectively. Figure~\ref{fig:st_deim_error} shows the error of the DEIM-augmented version with $p = 200$ along with that of the non-augmented version. There is almost no error sacrifice relative to the non-approximated method with this number of sample points. The timing for the method applied to the scalar transport problem is shown in Figure~\ref{fig:st_deim_timing}, along with the timing for a DEIM-augmented version of POD-Galerkin projection. Indeed, DEIM offers substantial speedup, both for the proposed method and for POD-Galerkin projection. 
\\

Two drawbacks of the DEIM-augmented version lead us to favor the non-augmented method in most cases. First, the DEIM-augmented version of the method relies on the initial condition and forcing being accurately approximated by sparse samplings, and these sparse samplings will only be accurate if the initial condition and forcing are similar in character to those in the training data. Second, the DEIM-augmented version is cumbersome to implement, requiring more precomputation of modes and matrices, relative to the non-augmented method.
\end{appendices}

\bibliographystyle{abbrv}
\bibliography{arXiv_template_bib}

\end{document}